\documentclass[aos,secthm,seceqn,nameyear,dvips]{arximspdf}
\usepackage{mathbh,dcolumn}
\usepackage{graphicx}

\doi{10.1214/10-AOS856}
\volume{39}
\issue{2}
\pubyear{2011}
\firstpage{803}
\lastpage{837}

\makeatletter
\newtheorem{lem}{Lemma}[section]
\newtheorem{prop}{Proposition}[section]
\newtheorem{cor}{Corollary}[section]

\newproclaim{rem}{Remark}[section]
\newproclaim{example}{Example}

\newcommand{\1}{\mathbh{1}}
\newcommand{\iint}{\int\!\!\!\int}

\newcolumntype{d}[1]{D{.}{.}{#1}}
\makeatother

\begin{document}
\begin{frontmatter}

\title{Estimation for L\'{e}vy processes from high frequency data
within a long time interval}
\runtitle{Estimation for L\'{e}vy processes}

\begin{aug}
\author[A]{\fnms{Fabienne} \snm{Comte}\corref{}\ead[label=e1]{fabienne.comte@parisdescartes.fr}}
and
\author[A]{\fnms{Valentine} \snm{Genon-Catalot}\ead[label=e2]{valentine.genon-catalot@parisdescartes.fr}}
\runauthor{F. Comte and V. Genon-Catalot}
\affiliation{University Paris Descartes, MAP5}
\address[A]{MAP 5 CNRS-UMR 8145\\
 Universit\'{e} Paris Descartes\\
45, rue des Saints-P\`{e}res\\
75006 Paris\\
France\\
\printead{e1}\\
\hphantom{E-mail:\ }\printead*{e2}} 
\end{aug}

\received{\smonth{4} \syear{2010}}
\revised{\smonth{9} \syear{2010}}

%
\begin{abstract}
In this paper, we study nonparametric estimation of the L\'{e}vy
density for L\'{e}vy processes, with and without Brownian component.
For this, we consider $n$ discrete time observations with step
$\Delta$. The asymptotic framework is: $n$ tends to infinity,
$\Delta=\Delta_n$ tends to zero while $n\Delta_n$ tends to infinity. We
use a Fourier approach to construct an adaptive nonparametric
estimator of the L\'{e}vy density and to provide a bound for the global
${\mathbb L}^2$-risk. Estimators of the drift and of the variance of
the Gaussian component are also studied. We discuss rates of
convergence and give examples and simulation results for processes
fitting in our framework.
\end{abstract}

%


\end{frontmatter}

\section{Introduction}
Let $(L_t, t \ge0)$ be a real-valued L\'{e}vy process, that is,
a~process with stationary independent increments and c\`{a}dl\`{a}g
sample paths. The distribution of $(L_t, t \ge0)$ is
completely specified by the characteristic function $\psi_{t}(u)=
{\mathbb E}(\exp{i u L_t})$ of the random variable $L_t$ which has the
form
\begin{equation} \label{fcg}
\qquad \psi_{t}(u)=\exp t\biggl(iu{\tilde b}-\frac12
u^2\sigma^2 +\int_{{\mathbb R}/\{0\}}\bigl(e^{iux} -1-iux1_{|x|\le1}\bigr)
N(dx)\biggr),
\end{equation}
where ${\tilde b}\in{\mathbb R}$, $\sigma^2\ge 0$ and $N(dx)$ is
a positive measure on ${\mathbb R}/\{0\}$ satisfying $\int_{{\mathbb
R}/\{0\}} x^2 \wedge1 N(dx) < \infty$ [see, e.g., Bertoin
(\citeyear{B1996}) or Sato (\citeyear{S1999})]. Thus, the statistical
problem for
L\'{e}vy processes is the estimation of its characteristic
triple $({\tilde b},\sigma^2, N)$ where appears a
finite-dimensional parameter $({\tilde b},\sigma^2)$ and an infinite-dimensional parameter $N$, the L\'{e}vy measure. In most recent
contributions, authors consider a discrete time observation of the
sample path, with regular sampling interval $\Delta$. Therefore,
statistical procedures are based on the i.i.d. sample
composed of the increments $(Z_k=Z_k^{\Delta}=
L_{k\Delta}-L_{(k-1)\Delta}, k=1, \ldots, n)$. In the general
case, the distribution of the r.v. $Z_k$ is not explicitly
given as a function of $({\tilde b},\sigma^2, N)$. This is why
authors rather use the relationship between the
characteristic function $\psi_{\Delta}$ of $Z_k$ and the
characteristic triple. Assuming that $N(dx)=n(x)\,dx$ admits a
density, several papers concentrate on the estimation of the
L\'{e}vy density under various assumptions on the characteristic
triple, including the case of ${\tilde b}=\sigma^2=0$ or assuming
stronger integrability conditions on the L\'{e}vy density [see, e.g.,
Watteel and Kulperger (\citeyear{WK2003}), Jongbloed and van der
Meulen (\citeyear{JvM2006}), van Es, Gugushvili and Spreij (\citeyear
{EGS2007}), Figueroa-L\'{o}pez (\citeyear{F2009})
and the references therein, Comte and Genon-Catalot
(\citeyear{CG2009}, \citeyear{CG2010a}, \citeyear{CG2010b})]. The joint
estimation of $({\tilde b},\sigma^2, N)$ is
investigated in Neumann and Reiss (\citeyear{NR2009}) or Gugushvili
(\citeyear{G2009}). The
methods and results differ according to the asymptotic point of
view. One may consider that the sampling interval $\Delta$ is
fixed and that $n$ tends to infinity (low frequency data). This
approach, which is quite natural, raises mathematical difficulties
and does not take into account the underlying continuous time
model properties. One may consider that $\Delta=\Delta_n$ tends
to $0$ as $n$ tends to infinity (high frequency data). Under
the assumption that $\Delta_n$ tends to $0$ within a fixed length
time interval ($n\Delta_n=t$ fixed), the estimation of $\sigma$
has been widely investigated for L\'{e}vy processes [see, e.g., Woerner
(\citeyear{W2006}), Barndorff-Nielsen, Shephard and Winkel (\citeyear{BSW2006}),
Jacod (\citeyear{J2007})]. However, the L\'{e}vy density cannot be identified
from observations within a~finite-length time interval. To
identify all parameters in the high-frequency context, one has to
assume both that $\Delta_n$ tends to $0$ and $n\Delta_n$
tends to infinity. This is the point of view adopted in this
paper. Our main focus is the nonparametric estimation of the
L\'{e}vy density $n(\cdot)$ by an adaptive deconvolution method which
generalizes the study of Comte and Genon-Catalot (\citeyear{CG2009}).
We also
study estimators of the other parameters. More precisely, we assume
that the L\'{e}vy density satisfies
{\renewcommand{\theequation}{H1}
\begin{equation}
\int_{\mathbb R} x^2 n(x) \,dx < \infty.
\end{equation}}

\vspace*{-\baselineskip}
\noindent For statistical purposes, this assumption, which was proposed in
Neumann and Reiss (\citeyear{NR2009}), has several useful consequences.
First, for all $t$, ${\mathbb E}L_t^2<+\infty$ and as $\int_{{\mathbb
R}}(e^{iux} -1-iux)n(x)\,dx$ is well defined, we get the following
expression for (\ref{fcg}):
\setcounter{equation}{1}
\begin{equation} \label{fc}
\psi_{t}(u)=\exp t\Biggl(iu b-\frac12
u^2\sigma^2 +\int_{{\mathbb R}}(e^{iux} -1-iux)
n(x)\,dx\Biggr),
\end{equation}
where $b={\mathbb E}L_1$ has a statistical meaning (contrary to
${\tilde b}$).

In Section \ref{prelim}, we present our main assumptions and some
preliminary properties. In Section \ref{sigmanul}, we assume that
$\sigma=0$ and study the estimation of the function $h(x)=x^2 n(x)$.
Using a sample of size $2n$, we build two collections of estimators
$(\hat h_m,\bar h_m)_{m>0}$ indexed by a cut-off parameter $m$. The
collections are obtained by Fourier inversion of two different
estimators of the Fourier transform $h^*$ of the function $h$. The
estimators of $h^*$ are built using empirical estimators of the
characteristic function $\psi_{\Delta}$ and its first two derivatives.
First, we give a bound for the ${\mathbb L}^2$-risk of $(\hat h_m,\bar
h_m)$ for fixed $m$. Then, introducing an adequate penalty,\vspace*{1pt} we propose
a data-driven choice of the cut-off parameter which yields an
estimator $(\hat h_{\hat m},\bar h_{\bar m})$ for each collection. The
${\mathbb L}^2$-risk of these estimators is studied. We discuss the
rates of convergence reached on Sobolev classes of regularity for the
function $h$. In Section~\ref{sigmanonnul}, we consider the general
case. To reach the L\'{e}vy density and get rid of the unknown
$\sigma^2$, we must now use derivatives of $\psi_{\Delta}$ up to the
order $3$ and we estimate the function $p(x)=x^3 n(x)$ developing the
Fourier inversion approach and adaptive choice of the cut-off parameter
as for $h$. It is worth stressing that the point of view of small
sampling interval is crucial to our study. Indeed, it helps obtaining
simple estimators of $\psi_{\Delta}$ and its successive derivatives
which are used to estimate the Fourier transform $p^*$ of $p$. Section
\ref{param} is devoted to the estimation of $(b, \sigma)$. We study
classical empirical means of the observations. This gives an estimator
of $b$ but cannot give estimators of $\sigma$. To estimate $\sigma$, we
consider power variation estimators, introduced in Woerner
(\citeyear{W2006}), Barndorff-Nielsen, Shephard and Winkel
(\citeyear{BSW2006}), Jacod (\citeyear{J2007}), A\"{i}t-Sahalia and
Jacod (\citeyear{AJ2007}),
under the asymptotic framework of high frequency data within a long
time interval. In Section \ref{exs}, we give examples of L\'{e}vy
models satisfying our set of assumptions. We provide numerical
simulation results in Section \ref{simu}. Section \ref{proofs} contains
the main proofs. In the \hyperref[Ap]{Appendix}, two classical results, used in
proofs, are recalled.



\section{Assumptions and preliminary properties} \label{prelim}

Let us consider the two functions
\[
h(x)=x^2 n(x),\qquad  p(x)= x^3 n(x),
\]
and the assumptions
{\renewcommand{\theequation}{H2}
\begin{equation}
(k)\hspace*{105pt}\int_{\mathbb R} |x|^k n(x) \,dx < \infty,\hspace*{105pt}
\end{equation}}

\[
\cases{
\mbox{(H3)}\quad  h \mbox{ belongs to } {\mathbb L}^{2}({\mathbb R})
\cr
\mbox{(H4)}\quad  \displaystyle \int x^8 n^2(x) \,dx =\int x^4 h^2(x) \,dx<\infty
}
\]
or
\[
\cases{
\mbox{(H5)}\quad  p \mbox{ belongs to } {\mathbb L}^{2}({\mathbb
R})\cr
\mbox{(H6)}\quad \displaystyle  \int x^{12} n^2(x) \,dx=\int x^6 p^2(x) \,dx <\infty.
}
\]
Assumption (H2)$(k)$ is a moment assumption. Indeed, according to
Sato [(\citeyear{S1999}), Section 5.25, Theorem 5.23], ${\mathbb E}|
L_t|^k<\infty$ is equivalent to\break $\int_{|x|>1} |x|^k n(x) \,dx<\infty$.
Below, for each stated result, the required value of $k$ is given.
Under (H1), the function $h$ is integrable and Section
\ref{sigmanul} is devoted to the nonparametric estimation of $h$
under the additional assumptions (H3)--(H4) when $\sigma^2=0$.
Assumption (H4) is only
required for the adaptive result. Under (H1)--(H2)$(3)$, the function
$p$ is integrable and Section \ref{sigmanonnul} concerns the
estimation of $p$ under (H5)--(H6) when \mbox{$\sigma^2\neq0$}.

Properties of the moments of $L_{\Delta}=Z_1^{\Delta}=Z_1$ for small
$\Delta$ are used in the proofs below.
\begin{lem}\label{ordremom2} Let $k\geq1$ be an integer and assume
\textup{(H1)--(H2)}$(k)$ with $k=3$ (or $k \ge3$). Then,
${\mathbb E}(|Z_1|^k)<+\infty$ and ${\mathbb E}(Z_1)=b\Delta$,
$\operatorname{Var}(Z_1)= \Delta(\sigma^2+\int x^2n(x)\,dx)$ and for $3\leq\ell\leq k$,
${\mathbb E}(Z_1^\ell) = \Delta c_\ell+ o(\Delta)$ where $c_\ell= \int
x^{\ell}n(x)\,dx$.
\end{lem}

Thus, under (H1), (H2)$(k)$, ${\mathbb E}(Z_1^{\ell}/\Delta)$ is
bounded for all $\ell\leq k$, for all $\Delta$.

In the sequel, results on the behavior of the characteristic function
$\psi_{\Delta}$ [see (\ref{fc})] for small $\Delta$ are needed.

\begin{lem}\label{bias1} Under \textup{(H1)}, $|\psi_{\Delta}(u)-1|\leq \Delta
|u|(c(u) + \sigma^2 |u|)$ where $c(u)=|b|+ |\int_0^u |h^*(v)|\,dv|$,
$h^*(v)= \int e^{ivx}h(x)\,dx$ denotes the Fourier transform of $h$. If
$h^*$ is integrable on ${\mathbb R}$, then
\[
|\psi_{\Delta}(u)-1|\leq \Delta
|u|(|b|+ |h^*|_1 + |u|\sigma^2) .
\]
\end{lem}
\begin{pf} By formula (\ref{fc}), under (H1), $\psi_{\Delta}$ is
$C^1$ with
$\psi'_{\Delta}(u)=\Delta\psi_\Delta(u) \times(\phi(u) -\sigma^2 u),$
where we have set, using that $e^{iux}-1= ix\int_0^u e^{ivx}\,dv$,
\setcounter{equation}{0}
\begin{equation}\label{phidu}
\phi(u)= i b -\int_0^u h^*(v)\, dv.
\end{equation}
We have $|\phi(u)|\le|b| +|\int_0^u|h^*(v)|\,dv|$
and by the Taylor formula, $\psi_\Delta(u) - 1= u\psi'_{\Delta}(c_u u)$
for some $c_u \in(0,1)$. The result follows.
\end{pf}

\section{Case of no Gaussian component} \label{sigmanul}
In this section, we consider the case $\sigma^2=0$ and focus on the
nonparametric estimation of $h$. For reasons that will appear below, we
suppose that we have at our disposal a
$2n$-sample,
$(Z_k)_{1\leq k\leq2n}$, with $Z_k=Z_k^{\Delta}=L_{k\Delta
}-L_{(k-1)\Delta}$. We assume that $\Delta=\Delta_n$ tends to $0$ and
$n\Delta_n$ tends to infinity. Hence, $\Delta$ and $Z_k$ depend on $n$.
However, to simplify notation, we omit the dependence on $n$ and
simply write $\Delta, Z_k$.

\subsection{Definition of estimators depending on a cut-off parameter}

For a~complex valued function $f$ belonging to ${\mathbb L}^1({\mathbb
R})$, we denote its Fourier transform by $f^*(u)=\int e^{iux}f(x)\,dx$.
For integrable and square integrable functions $f$, $f_1$, $f_2$, we use the
following notation:
\[
\|f\|=\int|f(x)|^2 \,dx, \qquad \langle f_1,f_2\rangle= \int f_1(x) {\bar f}_2(x)
\,dx
\]
(${\bar z}$ denotes the conjugate of the complex number $z$).
We have:
$
(f^*)^*(x)=2\pi f(-x)$ and $\langle f_1,f_2\rangle =1/(2\pi)\langle f^*_1,f^*_2\rangle .
$

By formula (\ref{fc}), under (H1), $\psi_{\Delta}$ is $C^2$ and we
have, as $\sigma^2=0$ [see (\ref{phidu})];
\[
\frac{\psi'_{\Delta}(u)}{\psi_\Delta(u)} = i\Delta\Biggl(b+ \int
\frac{e^{iux}-1}x h(x)\,dx\Biggr)= \Delta\phi(u).
\]
%
Derivating again gives
%
\begin{equation} \label{hstar}
h^*(u)= -\frac1\Delta \biggl(\frac{\psi''_{\Delta}(u)\psi_{\Delta}(u)-(\psi
'_\Delta(u))^2}{\psi_\Delta^2(u)}
\biggr),
\end{equation}
where, for all $u$, $\lim_{\Delta\rightarrow0} \psi_{\Delta}(u)=1$.
By splitting the $2n$-sample into two independent subsamples of $n$
observations, we introduce the following empirical unbiased estimators
of $\psi_{\Delta}, \psi'_{\Delta}, \psi''_{\Delta}$:
\[
\hat\psi_{\Delta,q}^{(j)}(u)= \frac1n \sum_{k=1+(q-1)n}^{qn}
(iZ_k)^j e^{iuZ_k},\qquad  j=0, 1, 2, q=1,2.
\]
We also define, based on the full sample, the estimator of $\psi
''_{\Delta}$
\[
\hat\psi^{(2)}_{\Delta}(u) = \frac1{2n}\sum_{k=1}^{2n} (iZ_k)^2 e^{iuZ_k}.
\]
We now build estimators of the Fourier transform $h^*$ of $h$.
Considering the expression
of $h^*$ in (\ref{hstar}), we replace $\psi_{\Delta}, \psi'_{\Delta},
\psi''_{\Delta}$ in the numerator by the empirical estimators built
on the two independent subsamples of size $n$. In the denominator, $\psi
_{\Delta}^2$ is simply replaced by $1$. This yields
\begin{equation}\label{hstarhat}
\hat h^*(u)= \frac1{\Delta} \bigl(\hat\psi^{(1)}_{\Delta,1}(u)\hat
\psi^{(1)}_{\Delta,2}(u)- \hat\psi^{(2)}_{\Delta,1}(u)\hat
\psi^{(0)}_{\Delta,2}(u )\bigr).
\end{equation}
Hence, using independence of the two subsamples,
\[
{\mathbb E}\hat h^*(u)= \frac1{\Delta} \bigl((\psi'_{\Delta}(u))^2- \psi
''_{\Delta}(u)
\psi_{\Delta}(u) \bigr)= h^{*}(u) + h^{*}(u)\bigl(\psi_{\Delta}^{2}(u)-1\bigr).
\]
Introducing a cut-off parameter $m$, we define an associated estimator
of $h$
\[
\hat h_m(x)=\frac1{2\pi} \int_{-\pi m}^{\pi m} e^{-iux}\hat
h^*(u)\,du.
\]
This means that $\hat h_m^*(u)=\hat h^*(u) 1_{[-\pi m, \pi m]}(u).$
By integration,
the following expression is available:
\[
\hat h_m(x)= \frac 1{ n^2 \Delta} \sum_{1\leq j,k\leq n}
(Z_k^2-Z_kZ_{n+j})\frac{\sin(\pi m(Z_k+Z_{j+n} -x))}{\pi(Z_k+Z_{j+n}-x)}.
\]

We also define another estimator of $h^*$ of $h$ by setting
\begin{equation}\label{hstarbar}
\bar h^*(u)=- \frac1{\Delta} \hat\psi^{(2)}_{\Delta}(u).
\end{equation}
Here, using (\ref{hstar}), we get
\begin{equation}\label{esphbar}
\quad {\mathbb E}\bar h^*(u)= -\frac1{\Delta} \psi''_{\Delta}(u)
= h^{*}(u) + h^{*}(u)\bigl(\psi_{\Delta}(u)-1\bigr)-\Delta\psi_{\Delta}(u)\phi
^2(u) .
\end{equation}
Thus, $\bar h^*$ is simpler but has an additional bias term. We set
\begin{equation}\label{hbar} \bar h_m(x)=\frac1{2\pi} \int_{-\pi
m}^{\pi m} e^{-iux}\bar
h^*(u)\,du= \frac{1}{2n \Delta}\sum_{k=1}^{2n}
Z_k^2 \frac{\sin(\pi m(Z_k -x))}{\pi(Z_k-x)} .
\end{equation}
%
\subsection{Risk for a fixed cut-off parameter}

Next, let us define
\[
h_m(x)=\frac1{2\pi}\int_{-\pi m}^{\pi m}
e^{-iux}h^*(u)\,du.
\]

Then we can prove the following result.
\begin{prop}\label{risk0} Assume that \textup{(H1)--(H2)}(4) and \textup{(H3)} hold. Then
\begin{eqnarray} \label{riskbound}
{\mathbb E}( \|\hat h_m-h\|^2)& \leq&\|
h_m-h\|^2 + 72{\mathbb E}(Z_1^4/\Delta)\frac m{n\Delta}\nonumber\\[-8pt]\\[-8pt]
&&{} +
\frac{4\Delta^2}{\pi} \int_{-\pi m}^{\pi m} u^2 c^2(u)|h^*(u)|^2\,du,\nonumber
\\
\label{riskbound2}
{\mathbb E}( \|\bar h_m-h\|^2) &\leq&
\|h_m-h\|^2 + {\mathbb E}(Z_1^4/\Delta)\frac m{n\Delta} \nonumber\\[-8pt]\\[-8pt]
&&{}+
\frac{2 \Delta^2}\pi \int_{-\pi m}^{\pi m}u^2 c^2(u)|h^*(u)|^2\,du +
C\Delta^2B_m,\nonumber
\end{eqnarray}
with $C$ a constant, $c(u)$ is defined in Lemma \ref{bias1}, $B_m
=(2/\pi)\int_{-\pi m}^{\pi m}|\phi(u)|^4\,du$ [see (\ref{phidu})]
satisfies $B_m=O(m)$ if $h^* \in{\mathbb L}_1({\mathbb R})$ and
$B_m=O(m^5)$ otherwise.
\end{prop}

\begin{rem} We stress that the estimator $\hat h_m$ is more complicated
to study, but $\bar h_m$ has an additional bias term.
\end{rem}
%
\subsection{Rates of convergence in Sobolev classes}
The following result concerns classes of functions $h$ belonging to
\begin{equation}\label{sobolev}
{\mathcal C}(a,L)= \Biggl\{ f\in({\mathbb L}^1\cap{\mathbb L}^2)({\mathbb
R}), \int(1+u^2)^a |f^*(u)|^2\,du \leq L\Biggr\}.
\end{equation}

\begin{prop}\label{rateofh} Assume that \textup{(H1)--(H2)}(4) and \textup{(H3)}
hold and that $h$ belongs to ${\mathcal C}(a,L)$ with $a>1/2$.
Consider the asymptotic setting where $n\rightarrow+\infty$,
$\Delta\rightarrow0$, $n\Delta\rightarrow+\infty$ and assume that
$m\le n\Delta$. If $n\Delta^{2}\leq1$, then, for the choice
$m=O((n\Delta)^{1/(2a+1)})$, we have
\[
{\mathbb E}( \|\hat
h_m-h\|^2) \leq O\bigl((n\Delta)^{-2a/(2a+1)}\bigr).
\]
If $a\geq1$, the
condition $n\Delta^2\leq1$ can be replaced by $n\Delta^3\leq1$.
The same result holds for ${\bar h}_m$.
\end{prop}

\begin{rem}
We can also discuss the case where $a \in(0, 1/2]$.
If $a \le1/2$, $|\int_0^u|h^*(v)|\,dv|=O(|u|^{1/2-a})$. Hence, the last
term in (\ref{riskbound}) is of order $\Delta^2m^{3-4a}$ which is less
than $m^{-2a}$ if $\Delta^2m^{3-2a}\leq1$ and thus $\Delta^2m^3\leq
1$. This requires $n\Delta^{5/3}\leq1$. The same holds for ${\bar h}_m$.
\end{rem}

Note that no lower bound result is available for this problem. A
benchmark for comparison could be the problem of density estimation for
i.i.d. observations without noise: if the density $f$ belongs to
${\mathcal C}(a,L)$, the optimal minimax rate is of order
$O(n^{-2a/(2a+1)})$ [see Ibragimov and Khas'minskij (\citeyear{IK1980})].

\subsection{Model selection}
The estimators ${\hat h}_m, {\bar h}_m$ are deconvolution estimators
that can also be described as minimum contrast estimators and
projection estimators. For details, the reader is referred to Comte and
Genon-Catalot (\citeyear{CG2009}, \citeyear{CG2010b}). For $m>0$, let
\[
S_m=\{f\in{\mathbb L}^2({\mathbb R}), \operatorname{support}(f^*) \subset
[-\pi m, \pi m]\}.
\]
The space $S_m$ is generated by an orthonormal basis, the sinus
cardinal basis, defined by
\[
\varphi_{m,j}(x)= \sqrt{m}\varphi(m x -j),\  j \in{\mathbb Z}, \qquad  \varphi
(x)=\frac{\sin{\pi x}}{\pi x}\ \bigl(\varphi(0)=1\bigr).
\]
This is due to the fact that
$
\varphi_{m,j}^{*}(u)= (e^{iuj/m}/\sqrt{m}) 1_{[-\pi m,\pi m]}(u), j\in
{\mathbb Z}.
$
For a function $f \in{\mathbb L}^2({\mathbb R})$, $f_m(x)= (2\pi)^{-1}
\int_{-\pi m}^{\pi m} e^{-iux}f^*(u)\,du$ is the orthogonal projection of
$f$ on $S_m$. Introducing, for a function $t \in S_m$,
\[
\gamma_n(t)= \|t\|^2 - \frac{1}{\pi}\langle{\hat h}^*, t^*\rangle= \|t\|^2 -
2\langle{\hat h}_m, t\rangle,
\]
we get
\[
{\hat h}_m= \arg\min_{t \in S_m}\gamma_n(t),
\]
and $\gamma_n({\hat h}_m)= -\|\hat h_m\|^2.$ We have
\[
{\hat h}_m= \sum_{j\in{\mathbb Z}}{\hat a}_{m,j} \varphi_{m,j}
\qquad \mbox{with } {\hat a}_{m,j}=\frac1{2\pi} \int_{-\pi
m}^{\pi m} {\hat h}^*(u)\varphi_{m,j}^{*}(-u)\,du
\]
and
$\|\hat h_m\|^2= 1/(2\pi) \int_{-\pi m}^{\pi m} |{\hat
h}^*(u)|^2 \,du $.
The coefficients ${\hat a}_{m,j} $ of the series as well as $\|\hat
h_m\|^2$ can be
explicitly computed by integration.
In the same way,~we set
\[
\Gamma_n(t)= \|t\|^2 - \frac{1}{\pi}\langle{\bar h}^*, t^*\rangle= \|t\|^2 -
2\langle{\bar h}_m, t\rangle,
\]
and obtain
\[
{\bar h}_m= \arg\min_{t \in S_m}\Gamma_n(t).
\]
Analogously, ${\bar h}_m$ has a series expansion on the sinus cardinal
basis with explicit coefficients and $\|\bar h_m\|^2$ has a closed-form
formula. We give the explicit expression of $\|\bar h_m\|^2$ which is
less cumbersome than $\|\hat h_m\|^2$:
\begin{equation}
\|\bar h_m\|^2= \frac{m}{4 n^2 \Delta^2} \sum_{1\le k,l\le2n} Z_k^2
Z_l^2 \varphi\bigl(m(Z_k - Z_l)\bigr).
\end{equation}
Now, we need to select the best $m$ as possible, in a set\vspace*{1pt} ${\mathcal
M}_n=\{m \in{\mathbb N},1\leq m\leq n\Delta\}=\{1, \dots, m_n\}$. For
the estimators ${\hat h}_m$,
we propose to take
\begin{equation}\label{hatm} \hat
m=\arg\min_{m\in{\mathcal M}_n} \bigl( -\|\hat h_m\|^2 + \operatorname{pen}(m) \bigr)
\end{equation}
with
\[
\operatorname{pen}(m)= \kappa\frac m{n\Delta^2} \Biggl(  \Biggl(\frac1n\sum_{k=1}^n
Z_k^2\Biggr)\biggl(\frac1n\sum_{k=n+1}^{2n} Z_k^2\biggr) +\frac1{n} \sum_{k=1}^n Z_k^4
\Biggr).
\]
The intuition for this choice is the following. The expression of
$\operatorname{pen}(m)$ is an estimator of the variance term of the risk
bound (\ref{riskbound}) as close as possible of the variance [see
(\ref{varhatstar})]. The term $-\|\hat h_m\|^2$ is an estimator of
$- \|h_m\|^2= \|h-h_m\|^2 - \|h\|^2$, which is up to a constant, the
bias term of the bound (\ref{riskbound}). This is why $\hat m$
mimics the optimal bias-variance compromise.

For the estimators ${\bar h}_m$,
we define
\begin{equation}\label{barm} \bar
m=\arg\min_{m\in{\mathcal M}_n} \Biggl( -\|\bar h_m\|^2 + \kappa' \frac
m{n\Delta^2} \Biggl( \frac1{2n} \sum_{k=1}^{2n} Z_k^4 \Biggr) \Biggr).
\end{equation}
The following result shows that the above data-driven choices of the
cut-off parameter lead to an automatic optimization of the risk.
\begin{thm}\label{main1}
Assume \textup{(H1)--(H2)(16)--(H3)--(H4)}. If, moreover, $h^*\in {\mathbb
L}^1({\mathbb R})$ and $n\Delta^3\leq1$, there exist  numerical
constants $\kappa,\kappa'$ such that
\begin{eqnarray*}
{\mathbb E}(\|\hat h_{\hat m}-h\|^2) &\leq& C\inf_{m\in{\mathcal
M}_n} \biggl(\|h-h_m\|^2 + \kappa\biggl(\Delta{\mathbb E}^2\biggl(\frac{Z_1^2}{\Delta
}\biggr)+{\mathbb E}\biggl(\frac{Z_1^4}\Delta\biggr)\biggr) \frac m{n\Delta} \biggr) \\
&&{} + \frac
{\Delta^2}\pi\int_{-\pi m_n}^{\pi m_n} u^2|h^*(u)|^2\,du +C\frac{\ln
^2(n\Delta)}{n\Delta},
\\
{\mathbb E}(\|\bar h_{\bar m}-h\|^2) &\leq& C\inf_{m\in{\mathcal
M}_n} \biggl(\|h-h_m\|^2 + \kappa'
{\mathbb E}\biggl(\frac{Z_1^4}\Delta\biggr) \frac m{n\Delta} \biggr) \\
&&{} + \frac{\Delta
^2}\pi\int_{-\pi m_n}^{\pi m_n} u^2|h^*(u)|^2\,du + \Delta^2 B_{m_n}
+C\frac{\ln^2(n\Delta)}{n\Delta},
\end{eqnarray*}
where $B_{m_n}=O(m_n)$ ($B_{m_n}$ is defined in Proposition \ref{risk0}).
\end{thm}

The numerical constants $\kappa, \kappa'$ have to be calibrated via
simulations [see discussion in Comte and Genon-Catalot (\citeyear{CG2009})].

By computations analogous to those in the proof of Proposition \ref
{rateofh}, we obtain the following corollary.
\begin{cor} Assume that the assumptions of Theorem \ref{main1} are
fulfilled. If, for some positive $L$, $h\in{\mathcal C}(a,L)$ with
$a>1/2$, then ${\mathbb E}(\|\hat h_{\hat m}-h\|^2)=O((n\Delta
)^{-2a/(2a+1)})$ provided that $n\Delta^2\leq1$. The same holds for
${\mathbb E}(\|\bar h_{\bar m}-h\|^2)$. If $a\geq1$, the constraint
$n\Delta^3\leq1$ is enough.
\end{cor}

\section{Study of the general case ($\sigma^2\neq0$)} \label{sigmanonnul}
In this section, we assume (H1)--(H2)($3$) and study the estimation of
the function
\[
p(x)=x^3n(x).
\]
We suppose that we have a sample of size $n$, $(Z_k)_{1\leq k\leq n}$,
$Z_k=L_{k\Delta}-L_{(k-1)\Delta}$.

\subsection{Definition of the estimators}
We compute the three first derivatives of $\psi_{\Delta}$ [see (\ref{phidu})]:
\[
\frac{\psi'_{\Delta}(u)}{\psi_\Delta(u)} = \Delta\Biggl(ib -u\sigma^2 + i\int
\frac{e^{iux}-1}x h(x)\,dx\Biggr)=\Delta\bigl(\phi(u)-u\sigma^2\bigr).
\]
Derivating again gives
\[
\frac{\psi''_{\Delta}(u)\psi_{\Delta}(u)-(\psi'_\Delta(u))^2}{(\psi
_\Delta(u))^2}= \Delta\bigl(\phi'(u)-\sigma^2\bigr)= -\Delta\Biggl(\sigma^2 + \int
e^{iux}x^2n(x)\,dx \Biggr),
\]
and last
%
\begin{eqnarray*}
p^*(u)
&=&\frac i\Delta \biggl(\frac{\psi_{\Delta}^{(3)}(u)}{\psi_{\Delta}(u)} - 3
\frac{\psi_{\Delta}''(u) \psi'_{\Delta}(u)}{ \psi_{\Delta}^2(u)} +
2\frac{[\psi'_{\Delta}(u)]^3}{\psi_{\Delta}^3(u)}
\biggr).
\end{eqnarray*}
%
Let
\[
\bar p^*(u)= \frac i{\Delta} \hat\psi_{\Delta}^{(3)}(u) \qquad \mbox{with }
\hat\psi_{\Delta}^{(3)}(u)= \frac1{n} \sum_{k=1}^{n}
(iZ_k)^3 e^{iuZ_k}.
\]
Then
\begin{equation} \label{pbar}
\bar p_m(x)=\frac1{2\pi} \int_{-\pi m}^{\pi m} e^{-iux}\bar p^*(u)\,du
=\frac{1}{n \Delta}\sum_{k=1}^{n} Z_k^3 \frac{\sin(\pi m(Z_k -x))}{\pi(Z_k-x)}.
\end{equation}
%
Let us set
\begin{equation}\label{phitilde}
{\tilde\phi}(u)= \phi(u) - u\sigma^2= ib - \int_0^u h^*(v)\,dv -u\sigma^2.
\end{equation}
Using $ \psi'_{\Delta}(u)= \Delta\psi_{\Delta}(u){\tilde\phi}(u)$ and
some computations, we get
\begin{eqnarray} \label{biaspbar}
\qquad {\mathbb E}\bar p^*(u)-p^*(u)
&=&\bigl(\psi_{\Delta}(u)-1\bigr) p^*(u) - 3 i \Delta\psi_{\Delta}(u){\tilde\phi
}(u) \bigl(\sigma^2+ h^*(u)\bigr) \nonumber\\[-8pt]\\[-8pt]
 &&{} +i\Delta^2 \psi_{\Delta
}(u)({\tilde\phi}(u))^3.\nonumber
\end{eqnarray}
%

\begin{rem} By a method analogous to the one used for $h$, considering
a sample of size $3n$, we can build another estimator of $p^*$ which is
less biased but more complicated to study.
\end{rem}

\subsection{Risk of the estimators}
The risk of the estimator with fixed cut-off parameter is bounded as follows.
\begin{prop}\label{risk2} Under \textup{(H1)--(H2)$(6)$} and \textup{(H5)},
%
\begin{eqnarray}\label{riskcomplet}
{\mathbb E}( \|\bar p_m-p\|^2) &\leq&
\|p-p_m\|^2 + {\mathbb E}(Z_1^6/\Delta)\frac m{n\Delta}
\nonumber\\[-8pt]\\[-8pt]
&&{} +
C \Biggl(\Delta^2\int_{-\pi m}^{\pi m} u^2(1+u^2)|p^*(u)|^2\,du + \Delta^2m^3 +
\Delta^4m^7 \Biggr),\nonumber
\end{eqnarray}
where $p_m(x)=(2\pi)^{-1}\int_{-\pi m}^{\pi m}
e^{-iux}p^*(u)\,du$
denotes the orthogonal projection of $p$ on $S_m$.
\end{prop}

We can state the result analogous to the one of Proposition \ref{rateofh}.
\begin{prop}\label{rateofp} Assume that \textup{(H1), (H2)($6$), (H5)} hold and
that $p$ belongs to ${\mathcal C}(a,L)$. Consider the asymptotic
setting where $n\rightarrow+\infty$,
$\Delta\rightarrow0$ and $n\Delta\rightarrow+\infty$.
If $n\Delta^{11/7} \leq1$, then
\[
{\mathbb E}( \|\bar p_{m}-p\|^2) \leq O\bigl((n\Delta)^{-2a/(2a+1)}\bigr).
\]
If $a\geq1/2$, the condition $n\Delta^{7/5}\leq1$ can be replaced by
$n\Delta^2\leq1$.
\end{prop}

\subsection{Model selection strategy}
The data driven selection of the best possible $m$ imposes here a
restricted collection of models. We choose
$ M_n=\{m \in{\mathbb N}/\{0\}, m\leq\sqrt{n\Delta}:={\mu}_n
\}$.

We can consider the estimator $ {\bar p}_{\bar m}$ where
\begin{eqnarray} \label{barmtilde} \bar
m=\arg\min_{m\in M_n} \bigl( -\|\bar p_m\|^2 + \overline{\operatorname{pen}}(m)
\bigr)\nonumber\\[-8pt]\\[-8pt]
\eqntext{\displaystyle \mbox{with } \overline{\operatorname{pen}}(m)= \kappa' \frac{m}{n\Delta^2}
\Biggl(\frac1{n} \sum_{k=1}^{n} Z_k^6 \Biggr).}
\end{eqnarray}

We can prove the following result.
\begin{thm}\label{main2}
Under assumptions \textup{(H1), (H2)($24$), (H5), (H6)} and with $n\Delta^2\leq
1$, there exists a numerical constant $\kappa$ such that (with $\mu
_n=\sqrt{n\Delta}$)
%
\begin{eqnarray*}
&&{\mathbb E}(\|\bar p_{\bar m}-p\|^2) \\
&&\qquad \leq C\inf_{m\in M_n} \biggl(\|
p-p_m\|^2 +\kappa'
{\mathbb E}\biggl(\frac{Z_1^6}{\Delta}\biggr)
\frac{m}{n\Delta} \biggr) \\
&&\qquad \quad {} + C \Biggl(\frac{\Delta^2}\pi\int_{-\pi\mu_n}^{\pi
\mu_n} u^2(1+u^2)|p^*(u)|^2\,du +\Delta^2 \mu_n^3 + \Delta^4 \mu_n^7+
\frac{\ln^2(n\Delta)}{n\Delta} \Biggr).
\end{eqnarray*}
\end{thm}


The consequence of Theorem \ref{main2} is that the adaptive estimators
reach automatically the expected rate of convergence when $p$ belongs
to a Sobolev class. This can be seen by computations analogous to those
of Proposition~\ref{rateofp}.

\section{Parameter estimation} \label{param}
Under (H1), the observed process may be written as $L_t= b t + \sigma
W_t + X_t$ where $(W_t)$ is a standard Brownian motion, $(X_t)$ is a L\'
{e}vy process, independent of $(W_t)$, of the form
\[
X_t=\int_{]0,t]}\int_{\mathbb R/\{0\}}x \bigl({\hat p}(ds,dx)- dsn(x)\,dx\bigr),
\]
where ${\hat p}(ds,dx)$ is the random jump measure of $(L_t)$ [and $(X_t)$].

If moreover
$
\int|x| n(x) \,dx <\infty,
$
then $L_t= b_0 t + \sigma W_t + \Gamma_t$ where $b_0= b- \int x n(x)
\,dx$ and
\[
\Gamma_t = \int_{]0,t]}\int_{\mathbb R}x {\hat p}(ds,dx)= X_t + t\int x
n(x)\, dx=\sum_{s\le t}\Gamma_s-\Gamma_{s_{-}}
\]
is of bounded variation on compact sets. We consider here a sample
of size $n$. By using empirical means of the data $Z_k^{\ell}$, it
is possible to obtain consistent and asymptotically Gaussian
estimators of $b$ ($\ell=1$) and, under suitable integrability
assumptions on the L\'{e}vy density, of $\int x^{\ell} n(x)\, dx$ for
$\ell\ge3$. But this method fails to estimate $\sigma$ for $\ell=2$
(see below). For this, one has to use another approach based
on power variations.

\subsection{Some small time properties}
To study estimators of $b$ and $\sigma$, small time properties of
moments of $L_\Delta$ are
needed. For simple moments, the result is stated in Lemma \ref
{ordremom2}. For absolute moments, we refer, for example, to
Figueroa-L\'{o}pez (\citeyear{F2008}): if $\int_{\{|x|>1\}}|x|^r
n(x)\,dx<+\infty$, and $r>2$, $\Delta^{-1}{\mathbb E}(|L_\Delta|^r)
\rightarrow\int|x|^r n(x)\,dx$ as $\Delta\rightarrow0$.
For the case of $|x|^r$ with $r<2$, we state the following proposition.
\begin{prop}\label{gamma}
\textup{(i)} Let $(\Gamma_t)$ be a L\'{e}vy process
with no continuous component and L\'{e}vy measure $n(\gamma)\,d\gamma$.
If $\int |\gamma| n(\gamma)\, d\gamma<\infty$, $b=\int\gamma n(\gamma)
\,d\gamma$ and for $r\le1$, $\int|\gamma|^r n(\gamma) \,d\gamma<\infty$.
There exists a constant $C$ such that, for all $\Delta$,
$ {\mathbb E}|\Gamma_{\Delta}|^r \le C \Delta.$
[Under the assumption, $(\Gamma_t)$ has finite mean and bounded
variation on compact sets.]
{\smallskipamount=0pt
\begin{longlist}[(iii)]
\item[(ii)] Let $X_t=B_{\Gamma_t}$ where $(\Gamma_t)$ is a subordinator with
L\'{e}vy
density $n_\Gamma$ satisfying $b=\int_0^{+\infty} \gamma n_\Gamma
(\gamma)
\,d\gamma<\infty$ and $(B_t)$ is a Brownian
motion independent of $(\Gamma_t)$.
The L\'{e}vy measure of $(X_t)$ has a
density given by
\begin{equation}\label{bgamma}
n_X(x)= \int_0^{+\infty} e^{-x^2/2\gamma}\frac{1}{\sqrt{2 \pi\gamma}}
n_\Gamma(\gamma)\,d\gamma.
\end{equation}
Consequently, if $C= \int_0^{+\infty} \gamma^{r/2} n_{\Gamma}(\gamma
)\,d\gamma<\infty$ with $r \le2$,
then
$
{\mathbb E}|X_{\Delta}|^r \le C \Delta.
$

\item[(iii)] Let $(X_t)$ be a L\'{e}vy process with no Gaussian component.
Then $X_{\Delta}/\break\sqrt{\Delta}$ converges to $0$ as $\Delta$ tends to
$0$ in probability and in ${\mathbb L}^{r}$ for all $r <2$.
\end{longlist}}
\end{prop}
%

\subsection{Estimator of $b$} Consider a L\'{e}vy process $(L_t)$
satisfying (H1) and set $Z_k= L_{k\Delta}-L_{(k-1)\Delta}$ as above.
Let us define the empirical means
\begin{equation}\label{estimofb}
\hat b = \frac{1}{n\Delta}\sum_{k=1}^n Z_k,\qquad  {\hat c}_{\ell}= \frac
{1}{n\Delta}\sum_{k=1}^n Z_k^{\ell}
\qquad \mbox{for } \ell\geq2.
\end{equation}
We prove now that $\hat b$, ${\hat c}_{\ell}, \ell\ge2$ are
consistent and asymptotically Gaussian estimators of the quantities
$b$, $c_{\ell}, \ell\ge2$ where
\[
c_2= \sigma^2 + \int x^2 n(x)\, dx,\qquad  c_{\ell}= \int x^{\ell} n(x) \,dx \qquad \mbox{for } \ell\ge3.
\]
\begin{prop} \label{emp} Assume \textup{(H1)} and $n$ tends to infinity, $\Delta
$ tends to $0$, $n\Delta$ tends to infinity.
\begin{longlist}[(ii)]
\item[(i)] Under \textup{(H2)}($2+\varepsilon$) for some positive $\varepsilon$,
\[
\sqrt{n\Delta}(\hat b - b ) \mbox{ converges in distribution to }
{\mathcal N}(0, c_2).\vspace{2pt}
\]
\item[(ii)] Under \textup{(H2)}($2(\ell+\varepsilon)$) for some positive $\varepsilon
$, and if $n\Delta^3$ tends to $0$,
$\sqrt{n\Delta}({\hat c}_{\ell}-c_{\ell})$ converges in distribution to
${\mathcal N}(0,c_{2\ell})$.
\end{longlist}
\end{prop}

We stress that this method provides an estimator of $b$ which is easy
to compute and very good in practice (see Section \ref{simu}),
but cannot provide an estimator of $\sigma^2$.\vspace{-2pt}

\subsection{Estimation of $\sigma$ with power variations}

Estimators of $\sigma$ based on power variations of $(L_t)$ have
been proposed and mostly studied in the case where $n\Delta=1$. They
are studied for high frequency data within a long time interval in
A\"{i}t-Sahalia and Jacod (\citeyear{AJ2007}). In the latter paper, the
context is
more general than ours, which implies that proofs are of high
complexity. For L\'{e}vy processes fitting in our set of assumptions,
we can derive the asymptotic properties of power variations
estimators with a specific proof given in Section \ref{proofs}.
Consider the family of estimators of $\sigma$ given by\vspace{-6pt}
\begin{equation}\label{sigmahatnr}
\hat\sigma(r) = \bigl[{\hat\sigma}_n^{(r)}\bigr]^{1/r} \qquad \mbox{with } {\hat\sigma
}_n^{(r)}= \frac{1}{m_r n \Delta^{r/2}} \sum_{k=1}^{n} |Z_k|^r,\vspace{-3pt}
\end{equation}
where $m_r={\mathbb E}|X|^r$ for $X$ a standard Gaussian variable
(recall that $Z_k= L_{k\Delta}-L_{(k-1)\Delta}$).

\begin{prop}\label{estimsigma}
As $n$ tends to infinity, $\Delta$ tends to 0 and $n\Delta$ tends~to
infinity, if $n\Delta^{2-r}=o(1)$, $\sqrt{n}({\hat
\sigma}_n^{(r)}-\sigma^r)$ converges in distribution to a ${\mathcal
N}(0,\break\sigma^{2r}(m_{2r}/m_r^2 -1))$
for:
\begin{longlist}[(ii)]
\item[(i)] $(L_t) $ a L\'{e}vy process satisfying \textup{(H1)} and such that
$\int|x| n(x) \,dx <\infty$ and $\int|x|^rn(x) \,dx <\infty$ for $r<1$.

\item[(ii)] $(L_t= bt +\sigma W_t +X_t) $, with $X_t=B_{\Gamma_t}$, where
$W,B, \Gamma$ are independent processes, $W,B$ are Brownian motions,
$\Gamma$ is a subordinator with L\'{e}vy measure $n_{\Gamma}$
satisfying $b=\int_0^{+\infty} \gamma n_\Gamma(\gamma) \,d\gamma<\infty
$ and $\int_0^{+\infty} \gamma^{r/2} n_\Gamma(\gamma) \,d\gamma<\infty$
for $r<1$.
\end{longlist}
Consequently, $\sqrt{n}(\hat\sigma(r)-\sigma)$ converges in
distribution to a ${\mathcal N}(0,(\sigma^2/r^2)(m_{2r}/\break m_r^2 -1))$.
\end{prop}

For other cases of L\'{e}vy processes, the result depends on the rate
of convergence to $0$ of ${\mathbb E}|X_{\Delta}|^{r}/\Delta^{r/2}$
[see Proposition \ref{gamma}(iii)] and will still hold if $\sqrt
{n\Delta}{\mathbb E}|X_{\Delta}|^{r}/\Delta^{r/2}$ tends to $0$.

\begin{rem}
It is worth noting that the rate of convergence is $\sqrt{n}$. For
$r=1$, the estimator ${\hat\sigma}_n^{(1)}$ is consistent but not
asymptotically Gaussian (because of its asymptotic bias). We have
implemented these estimators for $r=1/2$, $r=1/4$ (see Section \ref
{simu}) for processes satisfying $\int|x|^rn(x)\,dx<+\infty$ for all
positive $r$.
Note that we always give integrability conditions on ${\mathbb R}$ for
the L\'{e}vy density. This simplifies the presentation but induces some
redundancies. One should distinguish integrability conditions near $0$
and near infinity to avoid them.\vspace{-4pt}
\end{rem}

\section{Examples} \label{exs}
In this section, we give examples of models fitting in our framework.\vspace{-4pt}

\begin{example}\label{example1}
Drift${}+{}$Brownian motion${}+{}$Compound Poisson process.

Let\vspace{-10pt}
\begin{equation}\label{poisson} L_t= b_0t + \sigma W_t +
\sum_{i=1}^{N_t}Y_i,\vspace{-1pt}
\end{equation}
where $N_t$ is a Poisson process
with constant intensity $c$ and
$Y_i$ is a sequence of i.i.d. random variables with density $f$,
independent of the process
$(N_t)$. Then, $\sum_{i=1}^{N_t}Y_i$ is a compound Poisson process and
$(L_t)$ is a L\'{e}vy
process with L\'{e}vy density $n(x)=cf(x)$. Note that ${\mathbb
E}L_1=b=b_0+\int x n(x) \,dx$.
For the estimation of $p$, the rates that can be obtained depend on the
density $f$ provided
that $f$ satisfies the assumptions of Theorem \ref{main2}, which are
essentially here moment
assumptions for the r.v.'s $Y_i$. Any order can be obtained as shown in
Table \ref{fig0}
where rates are computed for $f$ a standard Gaussian, an exponential
with parameter
$1$ and a Beta distribution with parameters $(1,3)$ (for $p$ to be
regular enough).

\begin{table}
\tabcolsep=0pt
\caption{Rates for different ``Drift${}+{}$Brownian motion${}+{}$Compound Poisson processes''}
\label{fig0}\vspace{-5pt}
\begin{tabular*}{\tablewidth}{@{\extracolsep{\fill}}lccc@{}}
\hline
$\bolds{f(x)}$ & $\bolds{{\mathcal N}(0,1)}$ & $\bolds{{\mathcal E}(1)}$ & $\bolds{\beta(1,3)}$ \\
\hline
$p(x)=c x^3 f(x)$ & $\propto x^3e^{-x^2}$ & $\propto x^3e^{-x}\1_{x>0}$
& $\propto x^3(1-x)^2\1_{[0,1]}(x)$ \\
$p^*(u)$ & $\propto(u^3-3u)e^{-u^2/2}$ & $\propto1/(1-iu)^4$ &
$O(1/|u|^3)$ for large $|u|$\\
$\int_{|u|\geq\pi m} |p^*(u)|^2\, du $ & $O((\pi m)^5 e^{-(\pi m)^2})$ &
$O((\pi m)^{-7})$ & $O((\pi m)^{-5})$ \\
$\int_{|u|\leq\pi\mu_n}u^4|p^*(u)|^2\,du$ & $O(1)$ & $O(1)$ & $O(1)$ \\[3pt]
$\breve m$ (best choice of $m$) & $\sqrt{\log(n\Delta) - \frac52 \log
\log(n\Delta)}/\pi$ & $O((n\Delta)^{1/8})$ & $O((n\Delta)^{1/6})$\\[3pt]
Rate $\propto$ & $ \frac{\sqrt{\log(n\Delta)}}{n\Delta}$ & $(n\Delta
)^{-7/8}$ & $(n\Delta)^{-5/6}$ \\
\hline
\vspace{-20pt}
\end{tabular*}
\end{table}

As $\int|x|^r n(x) \,dx<\infty$ for all $ r<1$ (actually, for all $r\le
2$), estimation of
$\sigma$ is possible using $\hat\sigma(r)$ for any value of $0<r<1$ [provided
that $n\Delta^{2-r}=o(1)$].\vspace{-2pt}
\end{example}

\begin{example}\label{example2}
Drift${}+{}$Brownian motion${}+{}$L\'{e}vy--Gamma process.

Consider $L_t=b_0 t +\sigma W_t + \Gamma_t$ where
$(\Gamma_t)$ is a L\'{e}vy gamma process with parameters
$(\beta,\alpha)$, that is, is a subordinator such that, for all
$t>0$, $\Gamma_t$ has distribution Gamma with parameters $(\beta
t,\alpha)$ and density: $\alpha^{\beta t}x^{\beta t-1} e^{-\alpha
x}/\break\Gamma(\beta t) \1_{x\geq0}$. The L\'{e}vy density of $(L_t)$ is
$n(x)=\beta x^{-1}e^{-\alpha x}\1_{x>0}$. We have ${\mathbb
E}L_1=b=b_0+\int x n(x) \,dx$ and $p(x)=\beta x^2e^{-\alpha
x}\1_{x>0}$.

We find $p^*(u)= 2\beta/(\alpha- iu)^3$, $\int_{|u|\geq\pi m}
|p^*(u)|^2\,du=O(m^{-5})$ and\break $\int_{-\pi
\mu_n}^{\pi\mu_n} u^4|p^*(u)|^2\,du =O(1)$. Therefore, the rate for
estimating $p$ is $O((n\Delta)^{-5/6})$ for a choice $\breve
m=O((n\Delta)^{1/6})$.

As for all $r>0$, $\int x^r n(x) \,dx<\infty$, ${\hat\sigma}(r)$ is
authorized, for any value of $0< r<1$, to estimate $\sigma$.
\end{example}

\setcounter{example}{1}
\begin{example}[(Continued)]
Drift${}+{}$Brownian motion${}+{}$A specific class of subordinators.

Let $L_t=b_0 t +\sigma W_t + \Gamma_t$ where $(\Gamma_t)$ is a
subordinator of pure jump type with L\'{e}vy density of the form
$n(x)= \beta x^{\delta-1/2} x^{-1}e^{-\alpha x}\1_{x>0}$ with
$\delta>-1/2$ (thus, $\int xn(x)\,dx<\infty$). This class of
subordinators includes compound Poisson processes ($\delta>1/2$) and
L\'{e}vy Gamma processes ($\delta=1/2$). When $\delta>0$, the
function $xn(x)$ is both integrable and square integrable. This case
was discussed in Comte and Genon-Catalot (\citeyear
{CG2009}) where the estimation
of $xn(x)$, when $b_0=0$, $\sigma=0$, is studied. Here, we consider
the case $-1/2<\delta\le0$ which includes the L\'{e}vy Inverse
Gaussian process ($\delta=0$). Assumptions (H1)--(H6) are satisfied.
The function $p(x)=x^3n(x)$ can be estimated in presence (or not) of
additional drift and Brownian component. We can compute
\[
p^*(u)= \beta\frac{\Gamma(\delta+5/2)}{(\alpha-i u)^{\delta+5/2}}.
\]
Thus, $\int_{|u|\geq\pi m} |p^*(u)|^2 \,du =O(m^{-(2\delta+4)})$. As
$2\delta+1\le1$, $u^4|p^*(u)|^2$ is not integrable and we have $\Delta
^2 \int_{|u|\leq\pi\mu_n}u^4|p^*(u)|^2\,du= \Delta^2o(\mu_n)= o(\Delta
^{3/2})$. The best rate for estimating $p$ is $O((n\Delta)^{-(2\delta
+4)/(2\delta+5)})$ for a choice $\breve m=O((n\Delta)^{1/(2\delta
+5)})$. Note that $\Delta^{3/2}\le(n\Delta)^{-(2\delta+4)/(2\delta
+5)}$ for $n\Delta^2\le1$ and $-1/2<\delta\le0$.

We have $\int x^r n(x)\,dx<\infty$ for $r>1/2 - \delta$. Hence, to
estimate $\sigma$ using ${\hat\sigma}(r)$, we must choose $1/2 - \delta
<r<1$.
\end{example}

\begin{example}
Drift${}+{}$Brownian motion${}+{}$Pure jump
martingale.

Consider $L_t=bt + \sigma W_t + B_{\Gamma_t}$ where $W,B,\Gamma$ are
independent processes, $W,B$ are standard Brownian motion, and $\Gamma$
is a pure-jump subordinator with L\'{e}vy density $n_{\Gamma}(\gamma)=
\beta\gamma^{\delta-1/2} \gamma^{-1} e^{-\alpha\gamma}\1_{\gamma>0}$
as above (assuming $\delta>-1$). The L\'{e}vy density $n(\cdot)$ of $(L_t)$
[and of $(X_t=B_{\Gamma_t})$] is linked with $n_{\Gamma}$ [see (\ref
{bgamma})] and can be computed as the norming constant of a Generalized
Inverse Gaussian distribution
\[
n(x)=\frac{2\beta}{\sqrt{2\pi}} K_{\delta-1}\bigl(\sqrt{2\alpha}|x|\bigr)
\biggl(\frac{|x|}{\sqrt{2\alpha}}\biggr)^{\delta-1},
\]
where $K_{\nu}$ is a Bessel function of third kind (MacDonald
function) [see, e.g., Barndorff-Nielsen and Shephard (\citeyear{BS2001})].
For $\delta=1/2$, $B_{\Gamma_t}$ is a symmetric bilateral L\'{e}vy
Gamma process [see Madan and Seneta (\citeyear{MS1990}), K\"{u}chler
and Tappe
(\citeyear{KT2008})]. For $\delta=0$, $B_{\Gamma_t}$ is a normal inverse
Gaussian L\'{e}vy process [see Barndorff-Nielsen and Shephard
(\citeyear{BS2001})]. The relation (\ref{bgamma}) allows to
check that the function $p(x)=x^3 n(x)$ belongs to ${\mathbb L}^1
\cap{\mathbb L}^2$ and satisfies (H6) for $\delta>-3/4$. Moreover,
we can obtain
\[
p^*(u)= -i \beta \biggl(\frac{u^3\Gamma(\delta+5/2)}{(\alpha+
u^2/2)^{5/2}}-3 \frac{u\Gamma(\delta+3/2)}{(\alpha+u^2/2)^{3/2}} \biggr).
\]
Thus, $\int_{|u|\geq\pi m} |p^*(u)|^2 \,du =O(m^{-3})$ and $\Delta^2 \int
_{|u|\leq\pi\mu_n}u^4|p^*(u)|^2\,du=\break\Delta^2O(\mu_n)= O(\Delta
^{3/2})$. The best rate for estimating $p$ is $O((n\Delta)^{-3/4})$
obtained for $\breve m=O((n\Delta)^{1/4)})$. We have $\Delta^{3/2}\le
(n\Delta)^{-3/4}$ as $n\Delta^2\le1$.
As $\int\gamma^{r/2} n_{\Gamma}(\gamma)\,d\gamma<\infty$ for $r>1 -
\delta/2$, the estimation of $\sigma$ by ${\hat\sigma}(r)$ requires $1
- \delta/2<r<1$. Therefore, we must have $\delta>0$.
\end{example}

\section{Simulations}\label{simu}

In this section, we present numerical results for simulated L\'{e}vy
processes corresponding to Examples 1 and 2 (see Section
\ref{exs}). For these models, the functions $g(x)=xn(x)$, $h$ and
$p$ belong to ${\mathbb L}^1\cap{\mathbb L}^2({\mathbb R})$. Thus,
we can apply the method of Comte and Genon-Catalot (\citeyear{CG2009}), to
estimate $g$ when $b_0=0$, $\sigma=0$, and the method developed here
to estimate $h$ when $\sigma=0$ and $p$ when $\sigma\neq0$. We have
implemented the estimators $\bar h_{\bar m}$,
$\bar p_{\bar m}$ defined by (\ref{hbar})--(\ref{barm}) and
(\ref{pbar})--(\ref{barmtilde}). The numerical constant $\kappa'$
appearing in the penalties has been set to 7.5 for $g$, 4 for $h$
and 3 for $p$; its calibration is done by preliminary experiments.
The cutoff $\bar m$ is chosen among 100 equispaced values between 0
and 10.

\begin{figure}

\includegraphics{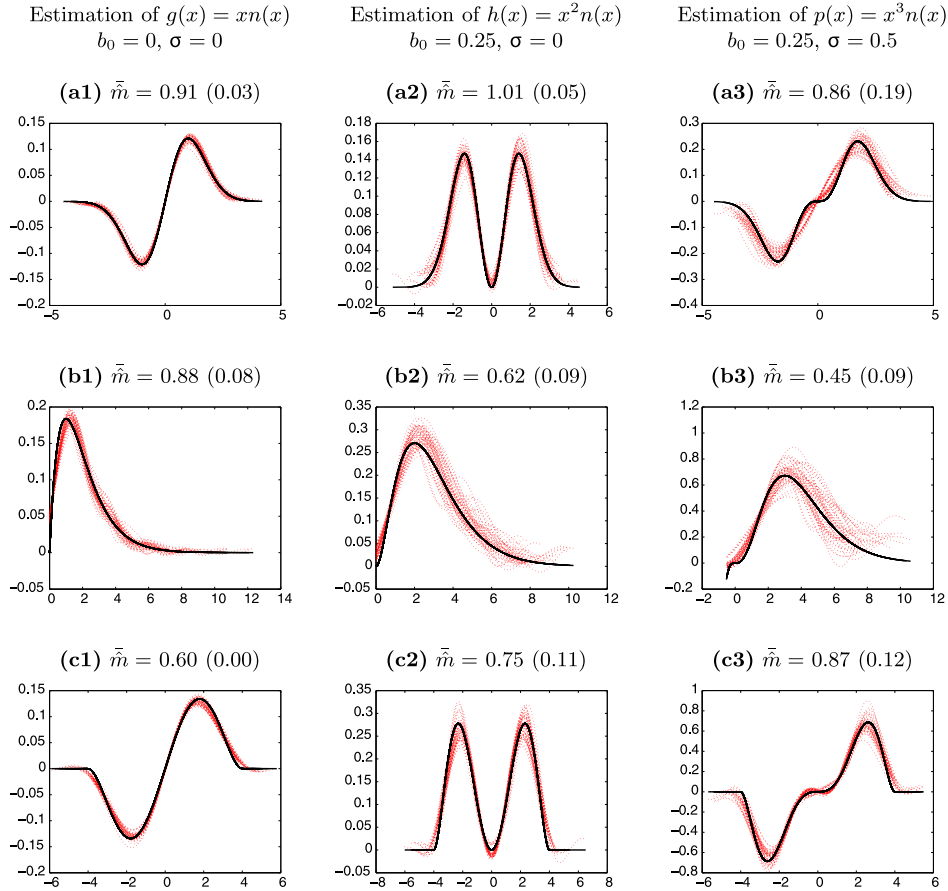}
\vspace{5pt}
\caption{Variability bands for the estimation of $g, h, p$ for a
compound Poisson process with Gaussian \textup{(first line)}, Exponential
${\mathcal E}(1)$ \textup{(second line)} and
$\beta(3,3)$ rescaled on $[-4,4]$ \textup{(third line)} $Y_i$'s, with $c=0.5$.
True (bold black line) and 50 estimated curves
(dotted red), $\Delta=0.05$, $n=5.10^4$.}\label{fig1}
\end{figure}

Figure \ref{fig1} shows estimated curves for models with jump part
coming from compound Poisson processes [see (\ref{poisson})] where
the $Y_i$'s are standard Gaussian, Exponential ${\mathcal E}(1)$, and
$\beta(3,3)$ rescaled on $[-4,4]$. The intensity
$c$ is equal to~0.5.

\begin{figure}

\includegraphics{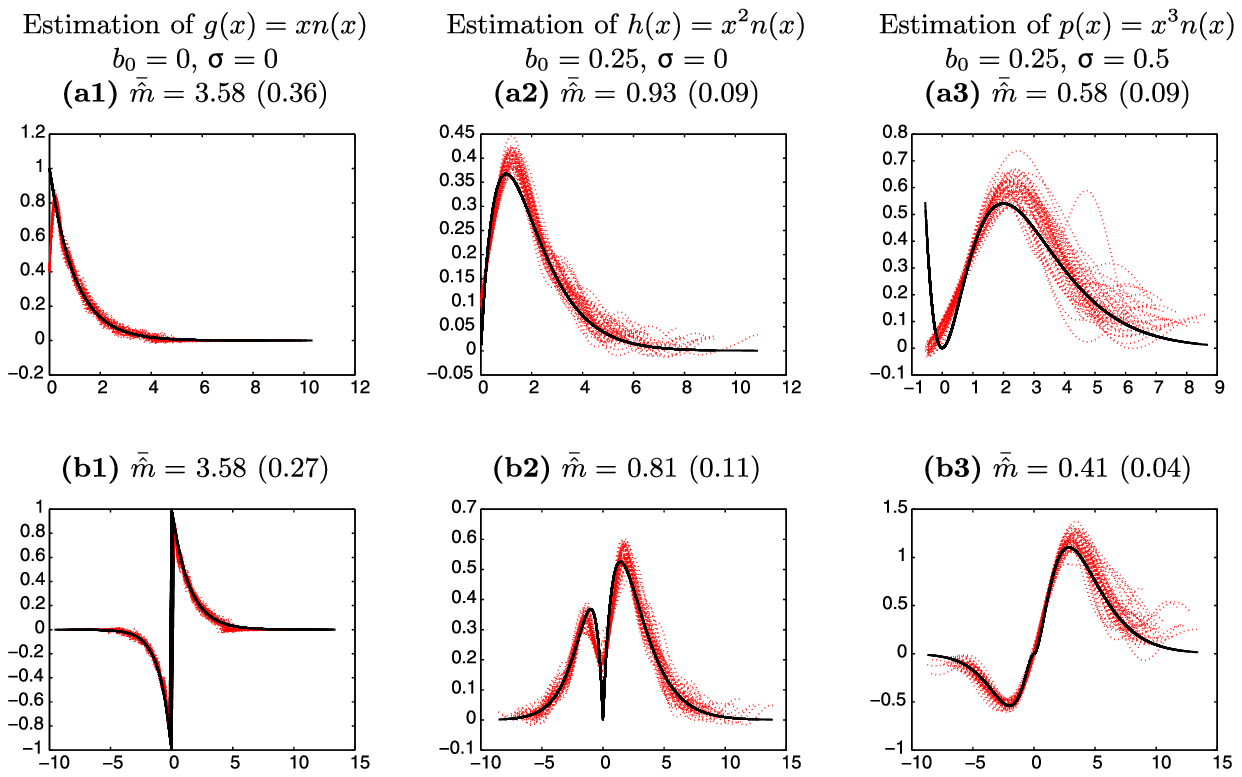}

\caption{Variability bands for the estimation of $g, h, p$ for jumps
from a
L\'{e}vy--Gamma process with $\beta=1, \alpha=1$ \textup{(first line)}, a
bilateral L\'{e}vy--Gamma process with
$(\beta,\alpha)=(0.7,1), (\beta',\alpha')=(1,1)$ \textup{(second line)}. True
(bold black line) and 50 estimated
curves (dotted red), $\Delta=0.05$, $n=5.10^4$.}\label{fig2}
\end{figure}

Figure \ref{fig2} shows estimated curves for jump part of L\'{e}vy
Gamma and bilateral L\'{e}vy Gamma type.
The bilateral L\'{e}vy Gamma process is the difference $\Gamma_t-\Gamma
'_t$ of two independent L\'{e}vy
Gamma processes.

On top of each graph, we give the mean value of the selected cutoff
with its standard deviation in parentheses. This value is
surprisingly small. As expected, the presence of a Gaussian
component deteriorates the estimation, which remains satisfactory on
the whole.

\begin{figure}

\includegraphics{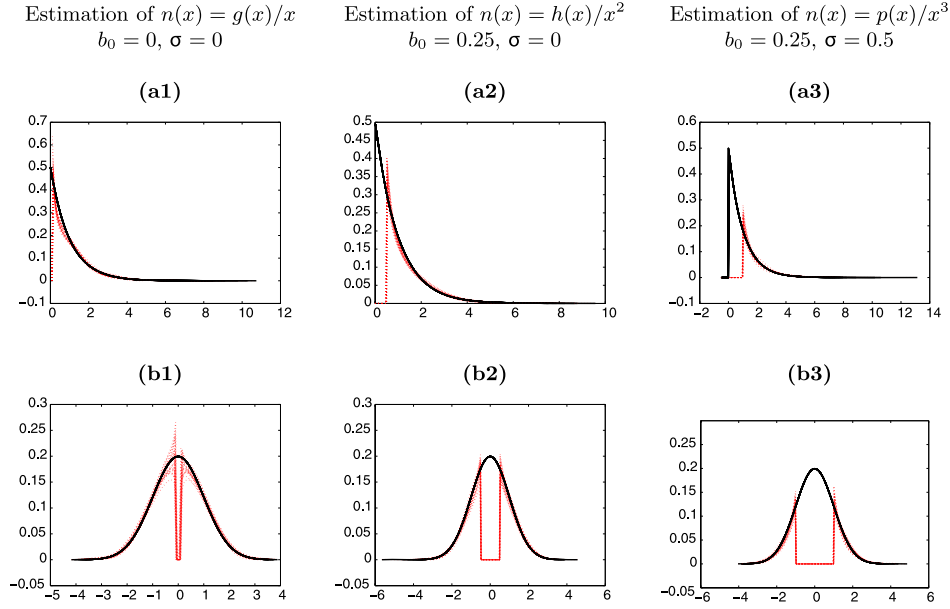}

\caption{Estimation of $n(\cdot)\1_{[-a,a]^c}$ with $a=0.1$ \textup{(first column)},
$a=0.5$ \textup{(second column)}, $a=1$ \textup{(third column)}. In all cases, $\lambda
=0.5$, $n=50\mbox{,}000$, $\Delta=0.05$; 25 estimated curves (thin dotted)${}+{}$the true (bold line).}
\label{estimn}\vspace*{10pt}
\end{figure}

We estimate the product of a power of $x$ and the L\'{e}vy density
whereas other authors
estimate $n(\cdot)$ on a compact set separated from the origin, see [12],
Figueroa-Lopez (\citeyear{F2009}).
Therefore, our point of view coincides with the usual one. Moreover we
have, an obvious
inequality; setting $\hat n(x)=\bar h(x)/x^2$ as $n(x)=h(x)/x^2$, we get
\[
{\mathbb E}\bigl(\|(\hat n-n)1_{{\mathbb R}/[-a,a]}\|^2\bigr) \leq\frac
1{a^2}{\mathbb E}(\|\bar
h-h\|^2).
\]
Analogous inequalities hold for $\hat n(x)=\hat g(x)/x$ or $\hat
n(x)=\bar p(x)/x^3$. In Figure~\ref{estimn}, we plot the estimator of
$n(\cdot)$ deduced
by dividing by the correct power of $x$ and by excluding an interval
$[-a,a]$ around zero. To obtain
correct representations, $a=0.1$ suits for $\hat g(x)/x$, $a=0.5$ for
$\bar h(x)/x^2$ and $a=1$
for $\bar p(x)/x^3$. The results are satisfactory and in accordance
with the difficulty of estimating $n(\cdot)$ without or with Gaussian
component.

Tables \ref{tab1} and \ref{tab2} show the means of the estimation
results for $b={\mathbb E}(L_1)= b_0+\int xn(x)\,dx$ [see (\ref
{estimofb})] and $\sigma$,
with standard deviations in parentheses.

The estimation of $b$ is good in all cases, and especially when $n\Delta
$ is large. The estimation of $\sigma$ is
clearly more difficult, with noticeable differences according to the
values of $n$ and $\Delta$. When $\Delta$ is not
small enough, the estimation can be heavily biased. In accordance with
the theory, when $r$ is
smaller, the estimator
of $\sigma$ is slightly better (smaller bias). Table \ref{tab3} shows
the values of $n\Delta^2$
and $n\Delta^{2-r}$, which should be small for the performance of the
estimator to be
satisfactory. It is worth noting that $\sigma$ is constantly over estimated.
\section{Proofs}\label{proofs}
\subsection{\texorpdfstring{Proof of Proposition \protect\ref{risk0}}{Proof of Proposition 3.1}}
 First, the Parseval
formula gives $\|\hat h_m - h\|^2= (1/(2\pi)) \|\hat h_m^*-h^*\|^2$
and we can note that $h^*(u)-h_m^*(u)=h^*(u)\1_{|u|\geq\pi m}$ is
orthogonal to $\hat h_m^*-h_m^*$ which has its support in $[-\pi m,
\pi m]$. Thus,
\[
\|\hat h_m - h\|^2= \frac1{2\pi} (\|h^*-h_m^*\|^2 + \|h_m^*-\hat
h_m^*\|^2).
\]
The first term $(1/(2\pi))
\|h^*-h_m^*\|^2=\|h-h_m\|^2$ is a classical squared bias term. Next,
\begin{eqnarray*}
\hat h_m^*(u)-h_m^*(u)&=& [\hat h_m^*(u)-{\mathbb
E}(\hat h_m^*(u)) ] + [{\mathbb E}(\hat h_m^*(u)) - h_m^*(u)] \\
&=&
[\hat h_m^*(u)-{\mathbb E}(\hat h_m^*(u)) ] + [\psi_{\Delta}^2(u)-1]
h^*(u)\1_{|u|\leq\pi m}.
\end{eqnarray*}
Bounding the norm of $\|\hat h_m^*-h_m^*\|^2$ by twice the sum
\begin{table}
\tabcolsep=0pt
\caption{Estimation of $(b,\sigma)$, $b_0=1$, the true value of $b$ in
parenthesis,\break
$\sigma=0.5$, $K=200$ replications}\label{tab1}
\begin{tabular*}{\tablewidth}{@{\extracolsep{\fill}}lcd{1.10}d{1.11}d{1.11}d{1.10}@{\hspace*{-3pt}}}
\hline
\textbf{Model} & $\bolds{(n,\Delta)}$ & \multicolumn{1}{c}{$\bolds{(5.10^4, 0.05)}$} &
\multicolumn{1}{c}{$\bolds{(5.10^4, 0.01)}$}
& \multicolumn{1}{c}{$\bolds{(5.10^4, 10^{-3})}$} & \multicolumn{1}{c@{\hspace*{-3pt}}}{$\bolds{(10^4,10^{-3})}$} \\
\hline
Poisson & $\hat b$ ($b=1$)& 1.000\ (0.02) & 0.997\ (0.04) & 0.995\ (0.123) & 1.001\ (0.280) \\
Gaussian & $\hat\sigma(1/2)$& 0.602\ (0.03) & 0.527\ (0.002) & 0.504
(0.002) & 0.504\ (0.005)\\
& $\hat\sigma(1/4)$ & 0.589\ (0.03) & 0.521\ (0.002) & 0.503\ (0.002) &
0.503\ (0.002) \\[3pt]
Poisson & $\hat b$ ($b=1.5$) & 1.502\ (0.05) & 1.502\ (0.051) & 1.494\
(0.142) & 1.461\ (0.359)\\
Exp(1) & $\hat\sigma(1/2)$& 0.611\ (0.003) & 0.530\ (0.003) & 0.505\
(0.002) & 0.505\
(0.005) \\
& $\hat\sigma(1/4)$ & 0.594\ (0.003) & 0.522\ (0.003) & 0.503\ (0.002)&
0.503\ (0.005) \\[3pt]
Gamma & $\hat b$ ($b=2$) & 2.001\ (0.02) & 2.000\ (0.05)& 1.998\
(0.177) & 2.018\ (0.335) \\
(1, 1) & $\hat\sigma(1/2)$& 0.705\ (0.004) & 0.562\ (0.003) & 0.512\
(0.002) & 0.513\ (0.005) \\
& $\hat\sigma(1/4)$ & 0.677\ (0.004) & 0.548\ (0.003) & 0.508\ (0.002) &
0.508\ (0.005) \\[3pt]
Bilateral & $\hat b$ ($b=1.4286$) & 1.426\ (0.035) & 1.4286\ (0.076) &
1.4493\ (0.264)& 1.405\ (0.619) \\
Gamma & $\hat\sigma(1/2)$& 0.862\ (0.005) & 0.628\ (0.004) & 0.526\
(0.003) & 0.526\ (0.006) \\
(0.7, 1), (1.1) & $\hat\sigma(1/4)$ & 0.798\ (0.004) & 0.593\ (0.003)
& 0.516\ (0.002) & 0.515\ (0.006) \\
\hline
\end{tabular*}\vspace*{7pt}
%
\caption{Estimation of $(b,\sigma)$, $b_0=1$, the true value of $b$ in
parenthesis,
$\sigma=1$, power variation method for estimation of $\sigma$, $K=200$
replications}\label{tab2}
\begin{tabular*}{\tablewidth}{@{\extracolsep{\fill}}lcd{1.10}d{1.10}d{1.10}d{1.10}@{\hspace*{-3pt}}}
\hline
\textbf{Model} & $\bolds{(n,\Delta)}$ & \multicolumn{1}{c}{$\bolds{(5.10^4, 0.05)}$} &
\multicolumn{1}{c}{$\bolds{(5.10^4, 0.01)}$}
& \multicolumn{1}{c}{$\bolds{(5.10^4, 10^{-3})}$} & \multicolumn{1}{c@{\hspace*{-3pt}}}{$\bolds{(10^4,10^{-3})}$} \\
\hline
Poisson & $\hat b$ (1)& 0.999\ (0.025) & 1.005\ (0.059) & 0.998\ (0.178) & 1.025\ (0.85) \\
Gaussian & $\hat\sigma(1/2)$& 1.082\ (0.005) & 1.026\ (0.004) & 1.006\
(0.004) & 1.005\ (0.009) \\
& $\hat\sigma(1/4)$ & 1.072\ (0.005) & 1.020\ (0.005)& 1.004\ (0.004) &
1.003\ (0.01) \\[3pt]
Poisson & $\hat b$\ (1.5) & 1.510\ (0.026) & 1.498\ (0.06) & 1.481\
(0.190) & 1.485\ (0.442) \\
Exp(1) & $\hat\sigma(1/2)$& 1.096\ (0.005) & 1.030\ (0.004) & 1.006\
(0.004) & 1.006\ (0.009) \\
& $\hat\sigma(1/4)$ &1.080\ (0.005) & 1.022\ (0.004) & 1.003\ (0.004) &
1.003\ (0.010) \\[3pt]
Gamma & $\hat b$\ (2) & 2.00\ (0.026) & 1.995\ (0.068) & 1.991\
(0.196) & 2.023\ (0.195) \\
(1, 1) & $\hat\sigma(1/2)$& 1.172\ (0.005) & 1.062\ (0.005) & 1.014\
(0.004) & 1.014\ (0.004) \\
& $\hat\sigma(1/4)$ & 1.152\ (0.005) & 1.050\ (0.005) & 1.010\ (0.005) &
1.010\ (0.004) \\[3pt]
Bilateral & $\hat b$\ (1.4286)& 1.425\ (0.04) & 1.431\ (0.10) &
1.429\ (0.28) & 1.492\ (0.63) \\
Gamma & $\hat\sigma(1/2)$& 1.330\ (0.006) & 1.136\ (0.005) & 1.033\
(0.005) & 1.033\ (0.01) \\
(0.7, 1),\ (1.1) & $\hat\sigma(1/4)$ & 1.284\ (0.006) & 1.105\ (0.005)
& 1.022\ (0.005)& 1.022\ (0.01) \\
\hline
\end{tabular*}\vspace*{7pt}
%
\caption{Values of $n,\Delta$, $n\Delta$, $n\Delta^2$, $n\Delta^{2-r}$
for $r=1/2$ and $r=1/4$}\label{tab3}
\begin{tabular}{@{}ld{4.0}d{3.0}d{1.2}d{1.2}@{}}
\hline
$\bolds{(n,\Delta)}$ & \multicolumn{1}{c}{$\bolds{(5.10^4, 0.05)}$} & \multicolumn{1}{c}{$\bolds{(5.10^4, 0.01)}$}
& \multicolumn{1}{c}{$\bolds{(5.10^4, 10^{-3})}$} & \multicolumn{1}{c@{}}{$\bolds{(10^4,10^{-3})}$} \\
\hline
$n\Delta$ & 2500 & 500 & 50 & 10 \\
$n\Delta^2$ & 125 & 5 & 0.05 & 0.01 \\
$n\Delta^{2-1/2}$ & 559 & 50 & 1.6 & 0.3\\
$n\Delta^{2-1/4}$ & 264 & 16 & 0.3 & 0.06 \\
\hline
\end{tabular}
\end{table}%
of the norms of the two elements of the decomposition, we get
\begin{eqnarray*}
{\mathbb E}(\|\hat h_m-h_m\|^2) &\leq& \frac1{\pi}
{\mathbb E} \Biggl(\int_{-\pi m}^{\pi m} |\hat h^*(u)-{\mathbb E}\hat
h^*(u)|^2 \,du \Biggr)\\
&&{}+ \frac1\pi\int_{-\pi m}^{\pi m} |\psi_{\Delta
}^2(u)-1|^2|h^*(u)|^2\,du
\\ &\leq& \frac1{\pi} \Biggl(\int_{-\pi m}^{\pi m} \operatorname{Var}(\hat h^*(u))\,du \Biggr)\\
&&{}+
\frac{4 \Delta^2}\pi\int_{-\pi m}^{\pi m} u^2c^2(u)|h^*(u)|^2\,du
\end{eqnarray*}
(see Lemma \ref{bias1} for the upper bound of $|\psi_{\Delta}(u)-1|$
and note that $|\psi_{\Delta}(u)|\le1$).
Now, we use the decomposition
\begin{eqnarray}\label{dechath}
&&\nonumber\Delta\bigl(\hat h^*(u)-{\mathbb E}(\hat h^*(u))\bigr)\\
\nonumber
&&\qquad =
\bigl(\hat\psi^{(1)}_{\Delta,1}(u)-\psi'_{\Delta}(u)\bigr) \bigl(\hat
\psi^{(1)}_{\Delta,2}(u)-\psi'_{\Delta}(u)\bigr)
\\
&&\qquad \quad {} + \bigl(\hat\psi^{(1)}_{\Delta,1}(u)-\psi'_{\Delta}(u)\bigr)\psi
'_{\Delta}(u) +
\bigl(\hat\psi^{(1)}_{\Delta,2}(u)-\psi'_{\Delta}(u)\bigr)\psi'_{\Delta}(u)\\
\nonumber
&&\qquad \quad  {}-\bigl( \hat\psi^{(2)}_{\Delta,1}(u)-\psi''_{\Delta}(u)\bigr)
\bigl(\hat\psi^{(0)}_{\Delta,2}(u)-\psi_{\Delta}(u)\bigr) \\
&&\qquad \quad {} -\bigl( \hat
\psi^{(2)}_{\Delta,1}(u)-\psi''_{\Delta}(u)\bigr)\psi_{\Delta}(u) - \bigl(\hat
\psi^{(0)}_{\Delta,2}(u)-\psi_{\Delta}(u)\bigr)\psi''_{\Delta}(u).\nonumber
\end{eqnarray}
Considering each term consecutively and exploiting the independence of
the samples, we obtain
\begin{eqnarray}  \label{varhatstar}
\operatorname{Var}(\hat h^*(u)) &\leq&\frac{6}{ \Delta^2} \biggl(\frac{{\mathbb
E}^2(Z_1^2)}{n^2} +
2\frac{{\mathbb E}^2(Z_1^2)}n + \frac{{\mathbb E}(Z_1^4)}{n^2} + 2\frac
{{\mathbb E}(Z_1^4)}{n} \biggr)\nonumber\\[-8pt]\\[-8pt]
&\leq&36
\frac{{\mathbb E}(Z_1^4/\Delta)}{n\Delta}.\nonumber
\end{eqnarray}
Thus, the first risk bound (\ref{riskbound}) is proved.
Analogously, we have
\begin{eqnarray*}
{\mathbb E}(\|\bar h_m-h\|^2) &\leq& \| h_m-h\|^2 + \frac1{\pi} \int
_{-\pi m}^{\pi m} |{\mathbb E}\bar h^*(u)- h^*(u)|^2 \,du\\
&&{}+ \frac1{\pi
}\int_{-\pi m}^{\pi m} \operatorname{Var}(\bar h^*(u))\,du.
\end{eqnarray*}
For the variance of $\bar h^*(u)$, we use:
$
\bar h^*(u)-{\mathbb E}\bar h^*(u) =
- \Delta^{-1}( \hat\psi_{\Delta}^{(2)}(u)-\psi''_{\Delta}(u)).
$
Thus,
\[
\operatorname{Var}(\bar h^*(u)) \le\frac{1}{2n \Delta} {\mathbb
E}(Z_1^{4}/\Delta).
\]
Next, for the bias of $\bar h^*(u)$, we use [see first (\ref{esphbar}) and then (\ref
{phidu})]
\[
|{\mathbb E}\bar h^*(u)- h^*(u)|^2 \le
2 |h^*(u)|^{2}| |\psi_{\Delta}(u)-1|^2 +2 \Delta^2 |\phi^4(u)| .
\]
Hence, there is an additional term in the risk bound equal to
\begin{equation}\label{bm} \frac{2}{\pi}\Delta^2 \int_{-\pi m}^{\pi m}
|\phi^4(u)|\,du=\Delta^2B_m.
\end{equation}
If $h^*$ is integrable, $|\phi(u)|\le C$  and $B_m=O(m)$. Otherwise,
$|\phi^4(u)|\le C|u|^4$  and $B_m=O(m^5)$.

\subsection{\texorpdfstring{Proof of Proposition \protect\ref{rateofh}}{Proof of Proposition 3.2}}
As $\|h-h_m\|^2=(1/\pi
) \int_{|u|\geq\pi m } |h^*(u)|^2\,du$, the definition of ${\mathcal
C}(a,L)$ implies clearly that $\|h-h_m\|^2\leq(L/2\pi) (\pi m)^{-2a}$.
The compromise between this term and the variance term of order
$m/(n\Delta)$ is standard: it leads to choose $m=O((n\Delta
)^{1/(2a+1)})$ and yields the order $O((n\Delta)^{-2a/(2a+1)})$.

For $a>1/2$, we have
\[
\biggl|\int_0^u|h^*(v)|\,dv\biggr|
\leq \sqrt{L\int(1+v^2)^{-a}\,dv}<+\infty.
\]
Therefore, $h^*$ is integrable and $|\phi(u)|\le|b|+|h^*|_1$.

The last term in the risk bound (\ref{riskbound}) is less than
\[
K\Delta^2\int_{-\pi m}^{\pi m} u^2 |h^*(u)|^2\,du \leq L\Delta^2(\pi
m)^{2(1-a)_+}.
\]

If $a\geq1$ and $n\Delta^3\le1$, we have $\Delta^2(\pi m)^{2(1-a)_+}=
\Delta^2\leq(n\Delta)^{-1}$.

If $a\in(1/2,1)$, the inequality $\Delta^2m^{2(1-a)} \leq m^{-2a}$ is
equivalent to $\Delta^2m^2\leq1$. As $m\leq n\Delta$, $\Delta^2m^2\leq
1$ holds if $n\Delta^{2}\leq1$.

For the additional bias term appearing in the risk bound of ${\bar
h}_m$, we have $B_m=O(m)$. Thus, $m\Delta^2\le m^{-2a}$ holds, for
$m=O((n\Delta)^{1/(2a+1)})$, if $m^{1+2a} \Delta^2=(n\Delta) \Delta^2
\le1$
which in turn holds if $n\Delta^{3}\le1$.

\subsection{\texorpdfstring{Proof of Theorem \protect\ref{main1}}{Proof of Theorem 3.1}}
We only study $\hat h_{\hat m}$ as the result for $\bar h_{\bar m}$ can
be proved analogously
(and is even simpler).

The proof is given in two steps. We define, for some $\varrho$,
$0<\varrho<1$,
\begin{eqnarray*}
\Omega_\varrho&:=& \biggl\{ \biggl| \frac{[(1/n\Delta)\sum_{k=1}^n Z_k^2][(1/n\Delta
)\sum_{k=n+1}^{2n} Z_k^2]}{({\mathbb E}(Z_1^2/\Delta))^2} -1 \biggr|\leq
\varrho/2 \biggr\}\\
&&{}\cap \biggl\{ \biggl| \frac{[(1/n\Delta)\sum_{k=1}^n
Z_k^4]}{({\mathbb E}(Z_1^4/\Delta))} -1 \biggr|\leq\varrho/2 \biggr\},
\end{eqnarray*}
so that ${\mathbb E}(\|\hat h_{\hat m}-h\|^2)={\mathbb E}(\|\hat
h_{\hat m}-h\|^2\1_{\Omega_\varrho}) +
{\mathbb E}(\|\hat h_{\hat m}-h\|^2\1_{\Omega_\varrho^c})$.

\textit{Step 1}. For the study of ${\mathbb E}(\|\hat h_{\hat
m}-h\|^2\1_{\Omega_\varrho^c})$, we refer to the analogous proof given
in Comte and Genon-Catalot (\citeyear{CG2009}) (see Section A4
therein). Using that ${\mathbb E}(Z_1^{16})<+\infty$, we can prove
${\mathbb E}(\|\hat h_{\hat m}-h\|^2\1_{\Omega_\varrho^c})\leq
C/(n\Delta)$. For this, we make use of the Rosenthal inequality [see
Hall and Heyde (\citeyear{HH1980})].


\textit{Step 2}. Study of ${\mathbb E}(\|\hat h_{\hat m}-h\|^2\1
_{\Omega_\varrho})$.

The proof relies on the following decomposition of $\gamma_n$:
\begin{eqnarray*}
\gamma_n(t)-\gamma_n(s)&=&\|t-h\|^2-\|s-h\|^2 +2\langle t-s, h\rangle
-\frac1\pi\langle{\hat h}^*,
t^*-s^*\rangle\\
&=& \|t-h\|^2-\|s-h\|^2 -2\nu_n(t-s)-2R_n(t-s),
\end{eqnarray*}
where
\[
\nu_n(t)=\frac{1}{2\pi} \langle\hat
h^*-{\mathbb E}(\hat h^*), t^*\rangle,\qquad  R_n(t)=\frac{1}{2\pi} \langle
{\mathbb E}(\hat h^*) - h^*, t^*\rangle.
\]
%
As $\gamma_n(\hat h_m)=-\|\hat h_m\|^2$, we deduce from (\ref{hatm})
that, for all
$m\in{\mathcal M}_n$,
\[
\gamma_n(\hat h_{\hat m})+ \operatorname{pen}(\hat m)\leq\gamma_n(h_m)+
\operatorname{pen}(m).
\]
This yields
\[
\|\hat h_{\hat m} - h\|^2 \leq\|h-h_m\|^2 + \operatorname{pen}(m) - \operatorname{pen}(\hat m)+2 \nu_n(\hat h_{\hat m} -h_m) + 2 R_n(\hat h_{\hat m}
-h_m).
\]
Then, for $\phi_n=\nu_n, R_n$, we use the inequality
\begin{eqnarray*}
2\phi_n(\hat h_{\hat m} -h_m) &\le& 2 \|\hat h_{\hat m} - h_m\| \sup
_{t\in S_m+S_{\hat m}, \|t\|=1} |\phi_n(t)|\\
&\le&\frac18 \| \hat h_{\hat m} -h_m\|^2 + 8 \sup_{t\in S_m+S_{\hat
m}, \|t\|=1} |\phi_n(t)|^2.
\end{eqnarray*}
Using that $\| \hat h_{\hat m} -h_m\|^2
\leq2 \| \hat h_{\hat m} -h\|^2 + 2 \| \hat h_m -h\|^2$ and some
algebra, we find
\begin{eqnarray}\label{risk}
\qquad  \frac14 \|\hat h_{\hat m} - h\|^2 & \leq&
\frac74 \|h-h_m\|^2 +
\operatorname{pen}(m) -\operatorname{pen}(\hat m)
\nonumber\\[-8pt]\\[-8pt]
&&{}
+ 8 \sup_{t\in S_m+S_{\hat m}, \|t\|=1} |R_n(t)|^2 + 8 \sup_{t\in
S_m+S_{\hat m},
\|t\|=1} |\nu_{n}(t)|^2.\nonumber
\end{eqnarray}
We have to study the terms containing a supremum, which are of
different nature. First, for $R_n(t)$, we have the following.

\begin{lem} $\!\!$We have:
$\sup_{t\in S_m+S_{\hat m}, \|t\|=1} |R_n(t)|^2 \le C\Delta^2\!\int_{-\pi
m_n}^{\pi m_n} u^2|h^*(u)|^2\,du.
$
\end{lem}

\begin{pf} We have $R_{n}(t)=\frac{1}{2\pi} \langle t^*, (1-\psi
_{\Delta}^2)h^* \rangle.$
By using Lemma \ref{bias1}, we find
\begin{eqnarray*}\nonumber\sup_{t\in S_m+S_{\hat m}, \|t\|=1} |
\langle t^*, (1-\psi_{\Delta}^2)h^* \rangle|^2
&\leq& \sup_{t\in S_{m_n}, \|t\|=1} | \langle t^*, (1-\psi_{\Delta
}^2)h^* \rangle|^2 \\
&\leq& 2\pi\bigl\|(1-\psi_\Delta^2)h^*\1_{[-\pi m_n, \pi m_n]}\bigr\|^2\\
&\leq&
C\Delta^2\int_{-\pi m_n}^{\pi m_n} u^2|h^*(u)|^2\,du.
\end{eqnarray*}
\upqed
\end{pf}

On the other hand, $\nu_n$ is decomposed: $\nu_n(t)=
\sum_{j=1}^4 \nu_{n,j}(t)+ r_n(t)$ with
\begin{eqnarray}\label{rn}
r_n(t)& = & \frac1{2\pi\Delta} \bigl\langle t^*,\bigl(\hat
\psi^{(1)}_{\Delta,1}(u)-\psi'_{\Delta}(u)\bigr) \bigl(\hat
\psi^{(1)}_{\Delta,2}(u)-\psi'_{\Delta}(u)\bigr)\bigr\rangle\nonumber\\[-8pt]\\[-8pt]
&&
{}-\frac1{2\pi\Delta} \bigl\langle t^*,\bigl( \hat
\psi^{(2)}_{\Delta,1}(u)-\psi''_{\Delta}(u)\bigr) \bigl(\hat
\psi^{(0)}_{\Delta,2}(u)-\psi_{\Delta}(u)\bigr) \bigr\rangle,\nonumber
\end{eqnarray}
and
\begin{eqnarray*}
\nu_{n,1}(t)& =& \frac1{2\pi\Delta} \bigl\langle t^*,\bigl( \psi''_{\Delta}-\hat
\psi^{(2)}_{\Delta,1}\bigr)\psi_{\Delta}\bigr\rangle,\qquad
\nu_{n,2}(t) = \frac1{2\pi\Delta} \bigl\langle t^*, \bigl(\psi_{\Delta} - \hat
\psi^{(0)}_{\Delta,2}\bigr)\psi_{\Delta}''\bigr\rangle.
\\
\nu_{n,3}(t) &=& \frac1{2\pi\Delta} \bigl\langle t^*, \bigl(\hat
\psi^{(1)}_{\Delta,1}-\psi'_{\Delta}\bigr)\psi'_{\Delta} \bigr\rangle,\qquad
  \nu_{n,4}(t) = \frac1{2\pi\Delta} \bigl\langle t^*, \bigl(\hat
\psi^{(1)}_{\Delta,2}-\psi'_{\Delta}\bigr)\psi'_{\Delta}\bigr\rangle.
\end{eqnarray*}

\begin{lem} We have:
${\mathbb E} (\sup_{t\in S_m+S_{\hat m}, \|t\|=1} |r_n(t)|^2 ) \le
\frac{C}{n}.
$
\end{lem}

\begin{pf} Using the independence of the subsamples, we can write
\begin{eqnarray}\label{borne2}
\nonumber
&&{\mathbb E} \Bigl(\sup_{t\in S_m+S_{\hat m}, \|
t\|=1} |r_n(t)|^2 \Bigr)\\
&&\qquad \leq
{\mathbb E} \Bigl(\sup_{t\in S_{m_n}, \|t\|=1} |r_n(t)|^2 \Bigr)\nonumber \\\nonumber
&&\qquad \leq
\frac1{2\pi^2 \Delta^2} {\mathbb E}\bigl[\bigl\| \bigl(\hat\psi^{(1)}_{\Delta,1}-\psi
'_{\Delta}\bigr) \bigl(\hat
\psi^{(1)}_{\Delta,2}-\psi'_{\Delta}\bigr)\1_{[-\pi m_n,\pi m_n]} \bigr\|^2 \\
&&\qquad \quad\hphantom{\frac1{2\pi^2 \Delta^2} {\mathbb E}\bigl[}
{} + \bigl\|\bigl(\hat
\psi^{(2)}_{\Delta,1}-\psi''_{\Delta}\bigr) \bigl(\hat
\psi^{(0)}_{\Delta,2}-\psi_{\Delta}\bigr)\1_{[-\pi m_n,\pi m_n]} \bigr\|^2\bigr]
\\ \nonumber
&&\qquad \leq \frac1{2\pi^2 \Delta^2} \int_{-\pi m_n}^{\pi m_n}
{\mathbb E}\bigl[\bigl|\hat\psi^{(1)}_{\Delta,1}(u)-\psi'_{\Delta}(u)\bigr|^2\bigr]
{\mathbb E} \bigl[\bigl|\hat
\psi^{(1)}_{\Delta,2}(u)-\psi'_{\Delta}(u)\bigr|^2\bigr] \,du \\
\nonumber
&&\qquad \quad
{}+ \frac1{2\pi^2 \Delta^2} \int_{-\pi m_n}^{\pi m_n} {\mathbb E}\bigl[\bigl|\hat
\psi^{(2)}_{\Delta,1}(u) -\psi''_{\Delta}(u)\bigr|^2\bigr] {\mathbb E}\bigl[\bigl|\hat
\psi^{(0)}_{\Delta,2}(u)-\psi_{\Delta}(u)\bigr|^2\bigr]\,du \\
&&\qquad \leq
 \frac{m_n}{\pi\Delta^2} \biggl( \frac{[{\mathbb E}(Z_1^2)]^2}{n^2} + \frac
{{\mathbb E}(Z_1^4)}{n^2} \biggr) \leq\frac{C}{n}\nonumber
\end{eqnarray}
because $m_n\leq n\Delta$ and ${\mathbb E}(Z_1^2)$ and ${\mathbb
E}(Z_1^4)$ have order $\Delta$.
\end{pf}

Now, the study of the $\nu_{n,j}$'s relies on Lemma \ref{Concent}. Let
us first study the process $\nu_{n,1}$.
We must split $Z_k^2=Z_k^2\1_{Z_k^2\leq k_n\sqrt{\Delta}} + Z_k^2\1
_{Z_k^2> k_n\sqrt{\Delta}}$ with $k_n$ to be defined later. This
implies that
$\nu_{n,1}(t) = \nu_{n,1}^P(t)+ \nu_{n,1}^R(t)$ ($P$ for Principal, $R$ for
residual) with
\begin{eqnarray} \label{nun1P}
\nu_{n,1}^P(t) = \frac1n
\sum_{k=1}^n [f_t(Z_k)-{\mathbb E}(f_t(Z_k))] \nonumber\\[-8pt]\\[-8pt]
\eqntext{\displaystyle \mbox{with }
f_t(z)=\frac1{2\pi\Delta} z^2\1_{z^2\leq k_n\sqrt{\Delta}} \langle
t^*, e^{iz\cdot}\psi_{\Delta}\rangle,}
\end{eqnarray}
and $\nu_{n,1}^R(t)=\nu_{n,1}(t) - \nu_{n,1}^P(t)$.
We prove the following results for $\nu_{n,1}$ and $\nu_{n,2}$.
\begin{prop}\label{concentP}
Under the assumptions of Theorem \ref{main1},
choose $k_n= C\frac{\sqrt{n}}{\ln(n\Delta)}$ and
\begin{equation}\label{p3mm} p(m,m')=4{\mathbb E}(Z_1^4/\Delta)\frac
{m\vee m'}{\Delta},
\end{equation}
then
\begin{eqnarray*}
&&{\mathbb E} \Bigl(\sup_{t\in S_m+S_{\hat m}, \|t\|=1} [\nu_{n,1}^{P}(t)]^2
-p(m,\hat m) \Bigr)_+ + {\mathbb E} \Bigl[\sup_{t\in S_{m_n}, \|t\|=1} \bigl|\nu
_{n,1}^{(R)}(t)\bigr|^2 \Bigr]\\
&&\qquad \leq C\frac{\ln^2(n\Delta)}{n\Delta},
\end{eqnarray*}
where $C$ is a constant.
\end{prop}
%
%
\begin{prop}\label{concentP4}
Under the assumptions of Theorem \ref{main1},
\[
{\mathbb E} \Bigl(\sup_{t\in S_m+S_{\hat m}, \|t\|=1} [\nu_{n,2}(t)]^2
-p(m,\hat m) \Bigr)_+ \leq\frac C{n\Delta},
\]
where $C$ is a constant.
\end{prop}

For both $\nu_{n,3}$ and $\nu_{n,4}$, which are similar, we have to
split again $Z_k=Z_k\1_{|Z_k|\leq k_n\sqrt{\Delta}} + Z_k\1_{|Z_k|>
k_n\sqrt{\Delta}}$ with the same $k_n$ as above.
We define $\nu_{n,j}(t) = \nu_{n,j}^P(t)+ \nu_{n,j}^R(t)$ as
previously, for $j=3,4$.
\begin{prop}\label{concentP12}
Under the assumptions of Theorem \ref{main1}, define for $j=3,4$
\begin{equation}\label{p12mm}
q(m,m')=4{\mathbb E}^2(Z_1^2/\Delta)\frac
{m\vee m'}{\Delta},
\end{equation}
then
\begin{eqnarray*}
&&{\mathbb E} \Bigl(\sup_{t\in S_m+S_{\hat m}, \|t\|=1} [\nu_{n,j}^{P}(t)]^2
-q(m,\hat m) \Bigr)_+ +{\mathbb E} \Bigl[\sup_{t\in S_{m_n}, \|t\|=1} \bigl|\nu
_{n,j}^{(R)}(t)\bigr|^2 \Bigr]\\
&&\qquad \leq C\frac{\ln^2(n\Delta)}{n\Delta},
\end{eqnarray*}
where $C$ is a constant.
\end{prop}

Now, on $\Omega_\varrho$, the following inequality holds (by bounding
the indicator by~1), for any choice of $\kappa$:
\[
(1-\varrho) \operatorname{pen}_{th}(m) \leq\operatorname{pen}(m) \leq(1+\varrho)\operatorname{pen}_{th}(m),
\]
where $\operatorname{pen}_{th}(m)={\mathbb E}(\operatorname{pen}(m))$. It follows from
(\ref{risk}) that
\begin{eqnarray}\label{t2}
\frac14 {\mathbb E}( \|\hat h_{\hat m} - h\|
^2 \1_{\Omega_\varrho}) & \leq& \frac74 \|h-h_m\|^2 + \operatorname{pen}_{th}(m) -{\mathbb E}(\operatorname{pen}(\hat m)\1_{\Omega_\varrho})
\nonumber\\
&&{}
+ C\Delta^2\int_{-\pi m_n}^{\pi m_n} u^2|h^*(u)|^2\,du \\
&&{}+ 8 {\mathbb
E}\Bigl(\sup_{t\in S_m+S_{\hat m},
\|t\|=1} |\nu_{n}(t)|^2 \1_{\Omega_\varrho}\Bigr).\nonumber
\end{eqnarray}
Recalling that
\[
\nu_{n}(t)=r_n(t) + \nu_{n,1}^P(t) + \nu_{n,1}^R(t) + \nu_{n,2}(t)+ \nu
_{n,3}^P(t)+ \nu_{n,3}^R(t)+ \nu_{n,4}^P(t)+ \nu_{n,4}^R(t),
\]
we have
\begin{eqnarray} \label{t1}
\nonumber
&&{\mathbb E}\Bigl(\sup_{t\in S_m+S_{\hat m},
\|t\|=1} |\nu_{n}(t)|^2 \1_{\Omega_\varrho}\Bigr)\\
&&\qquad \leq 8 \biggl( \frac C{n\Delta
} + \sum_{j\in\{1,3,4\}} {\mathbb E}\Bigl(\sup_{t\in S_m+S_{\hat m}, \|t\|
=1} |\nu_{n,j}^P(t)|^2 \1_{\Omega_\varrho}\Bigr)\nonumber\\[-8pt]\\[-8pt]
&&\qquad \quad\hspace*{65pt}
{}  +{\mathbb E}\Bigl(\sup_{t\in S_m+S_{\hat m}, \|t\|=1} |\nu_{n,2}(t)|^2 \1_{\Omega
_\varrho}\Bigr) \biggr) \nonumber \\
&&\qquad \leq
8\biggl( \frac{C'}{n\Delta} + 2{\mathbb E}\bigl[\bigl(p(m,\hat m) + q(m,\hat m)\bigr) \1
_{\Omega_\varrho}\bigr] \biggr).\nonumber
\end{eqnarray}
We note that $p(m,m')+q(m,m')= \frac1{4\kappa} (\operatorname{pen}_{th}(m) +
\operatorname{pen}_{th}(m')).$
Thus,
\begin{eqnarray*}
&&\operatorname{pen}_{th}(m) -{\mathbb E}(\operatorname{pen}(\hat m)\1_{\Omega_\varrho})+
128 {\mathbb E}\bigl[\bigl(p(m,\hat m) + q(m,\hat m)\bigr) \1_{\Omega_\varrho}\bigr] \\
&&\qquad \leq \operatorname{pen}_{th}(m) -(1-\varrho) {\mathbb E}(\operatorname{pen}_{th}(\hat
m)\1_{\Omega_\varrho})+ \frac{32}{\kappa} {\mathbb E}\bigl[\bigl(\operatorname{pen}_{th}(m) + \operatorname{pen}_{th}(\hat m)\bigr) \1_{\Omega_\varrho}\bigr] \\
&&\qquad \leq \biggl(1+ \frac{32}{\kappa}\biggr) \operatorname{pen}_{th}(m) + \biggl(\frac{32}{\kappa
}-(1-\varrho)\biggr) {\mathbb E}[\operatorname{pen}_{th}(\hat m) \1_{\Omega_\varrho}].
\end{eqnarray*}
Therefore, we choose $\kappa$ such that $(32/\kappa-(1-\varrho)) \leq
0$, that is $\kappa\geq32/(1-\varrho)$. This together with (\ref{t2})
and (\ref{t1}) yields
\begin{eqnarray*} \frac14 {\mathbb E}( \|\hat h_{\hat m} -
h\|^2 \1_{\Omega_\varrho}) & \leq& \frac74 \|h-h_m\|^2 + (2-\varrho)
\operatorname{pen}_{th}(m) \\ && {}+ C\Delta^2\int_{-\pi m_n}^{\pi m_n}
u^2|h^*(u)|^2\,du + \frac{C''}{n\Delta}.
\end{eqnarray*}

\subsection{\texorpdfstring{Proof of Propositions \protect\ref{concentP}--\protect\ref{concentP12}}%
{Proof of Propositions 8.1--8.3}}
\mbox{}

\begin{pf*}{Proof of Proposition \ref{concentP}}
Let $m''=m\vee m'$, and note that $S_m+S_{m'}=S_{m''}$. We evaluate the constants
$M,H,v$ to apply Lemma \ref{Concent} to $\nu_{n,1}^{P}(t)$ [see
(\ref{nun1P})]:
\begin{eqnarray*}
\sup_{z\in{\mathbb R}}|f_t(z)|&\leq& \frac
{k_n}{2\pi\sqrt{\Delta}} \sup_{z\in{\mathbb R}} \biggl|\int_{-\pi m''}^{\pi
m''} t^*(-u)e^{iuz}\psi_{\Delta}(u)\,du\biggr|\\
&\leq& \frac{k_n}{2\pi\sqrt{\Delta}}\int_{-\pi m''}^{\pi m''}
|t^*(u)|\,du \leq
\frac{k_n}{2\pi\sqrt{\Delta}}\Biggl(2\pi m'' \int_{-\pi m''}^{\pi m''}
|t^*(u)|^2\,du\Biggr)^{1/2} \\ &= &
\frac{k_n}{\sqrt{\Delta}}(m'')^{1/2} \|t\|= \frac{k_n\sqrt{m''}}{\sqrt
{\Delta}}:=M.
\end{eqnarray*}
Moreover,
\begin{eqnarray*}
{\mathbb E}\Bigl(\sup_{t\in S_m+S_{m'}, \|t\|=1} \bigl[\nu_{n,1}^{P}(t)\bigr]^2\Bigr)
&\leq&
\frac1{2\pi n\Delta^2} \int_{-\pi m''}^{\pi m''} {\mathbb E}(Z_1^4)\psi
_{\Delta}^2(u)\,du \\
&\leq&\frac{m'' {\mathbb E}(Z_1^4/\Delta)}{n\Delta}:=H^2.
\end{eqnarray*}
The most delicate term is $v$:
\begin{eqnarray*}
\operatorname{Var}(f_t(Z_1)) &=& \frac1{4\pi^2\Delta^2}
{\mathbb E} \biggl( Z_1^4
\1_{Z_1^2\leq k_n\sqrt{\Delta}} \biggl|\iint e^{ixZ_1}t^*(-x)\psi_{\Delta
}(x)\,dx\biggr|^2 \biggr)\\
&\leq& \frac1{4\pi^2\Delta^2}
{\mathbb E} \biggl(Z_1^4 \iint
e^{i(x-y)Z_1}t^*(-x)t^*(y)\psi_{\Delta}(x)\psi_{\Delta}(-y)\,dx\,dy \biggr)\\
&=& \frac1{4\pi^2\Delta^2}\iint\psi^{(4)}_{\Delta}(x-y)
t^*(-x)t^*(y)\psi_{\Delta}(x)\psi_{\Delta}(-y)\,dx\,dy,
\end{eqnarray*}
where we recall that $\psi^{(4)}_{\Delta}(x)={\mathbb E}( Z_1^4 e^{ixZ_1}).$
Making use of the basis $(\varphi_{m'',j}, j\in{\mathbb Z})$ of
$S_{m''}$, we have $t=\sum_{j\in{\mathbb Z}} t_j \varphi_{m'',j}$ with
$\|t\|^2=\sum_{j\in{\mathbb Z}}t_j^2=1$,
\begin{eqnarray}\label{presquev}\nonumber
\hspace*{30pt}\operatorname{Var}(f_t(Z_1)) &\leq&
\frac1{4\pi^2\Delta^2}\sum_{j,k\in{\mathbb Z}} t_{j}t_{k}\iint\psi
^{(4)}_{\Delta}(x-y)
\varphi_{m'',j}^*(-x)\varphi_{m'',k}^*(y)\\
\nonumber
&&\hphantom{\frac1{4\pi^2\Delta^2}\sum_{j,k\in{\mathbb Z}} t_{j}t_{k}\iint}
{}\times\psi_{\Delta}(x)\psi_{\Delta
}(-y) \,dx\,dy
\\\nonumber&\leq& \frac1{4\pi^2\Delta^2} \biggl(\sum_{j,k \in{\mathbb
Z}} \biggl|\iint\psi^{(4)}_{\Delta}(x-y) \varphi_{m'',j}^*(-x)\varphi
_{m'',k}^*(y)\\
\nonumber
&&\hspace*{114pt}
{}\times\psi_{\Delta}(x)\psi_{\Delta}(-y)\,dx\,dy \biggr|^2 \biggr)^{1/2}
\\\nonumber&= & \frac1{4\pi^2\Delta^2} \biggl(\iint_{[-\pi m'', \pi
m'']^2} \bigl|\psi^{(4)}_{\Delta}(x-y)\bigr|^2
|\psi_{\Delta}(x)|^2\\
\nonumber
&&\hspace*{120pt}
{}\times|\psi_{\Delta}(-y)|^2 \,dx\,dy \biggr)^{1/2},
\\
\hspace*{30pt}\operatorname{Var}(f_t(Z_1)) &\leq & \frac1{4\pi^2\Delta^2} \biggl(\iint
_{[-\pi m'', \pi m'']^2} \bigr|\psi^{(4)}_{\Delta}(x-y)\bigr|^2
\,dx\,dy \biggr)^{1/2}\nonumber\\[-8pt]\\[-8pt]
 &\leq & \frac{\sqrt{2\pi m''}}{4\pi^2\Delta^2} \biggl(\int
_{[-2\pi m'', 2\pi m'']} \bigl|\psi^{(4)}_{\Delta}(z)\bigr|^2 \,dz \biggr)^{1/2}.\nonumber
\end{eqnarray}
Therefore, we need to study $\int_{[-2\pi m'', 2\pi m'']} |\psi
^{(4)}_{\Delta}(z)|^2 \,dz$.
Recall that $\phi(u)= ib -\int_0^u h^*(v)\,dv$. We have
\[
\psi^{(4)}_{\Delta}=\Delta\bigl[\phi^{(3)}+\Delta\bigl(4\phi\phi''+3(\phi')^2\bigr)
+6\Delta^2\phi'\phi^2+\Delta^3\phi^4\bigr]
\psi_{\Delta},
\]
where
\begin{eqnarray*}
\phi'(u)&=&-h^*(u),\qquad  \phi''(u)=-i \int e^{iux} x^3 n(x)\,dx,\\
\phi^{(3)}(u)&=&
\int e^{iux} x^4 n(x)\,dx
\end{eqnarray*}
satisfy: $\int|\phi'(u)|^2 \,du = \|h\|^2$, $|\phi'(u)|\le|h|_1$  and
thanks to (H4), the
Parseval equality yields
\begin{eqnarray*}
\int|\phi''(u)|^2 \,du&=& \int x^6 n^2(x)\,dx=\int x^2 h^2(x)\,dx,
\\
\int|\phi^{(3)}(u)|^2 \,du&=&\int x^8 n^2(x) \,dx=\int x^4 h^2(x)\,dx.
\end{eqnarray*}
By assumption, $h^*$ is in ${\mathbb L}_1({\mathbb R})$, thus, $|\phi
(u)|\le|b|+|h^*|_1:=M_\phi$.
Therefore,
\[
|\psi^{(4)}_{\Delta}|^2\leq C \Delta^2\bigl(\bigl|\phi^{(3)}\bigr|^2 + \Delta^2 \bigl((\phi
'')^2 + (\phi')^4\bigr) + \Delta^4 (\phi')^2+\Delta^6\bigr),
\]
where $C$ is a constant depending on $M_\phi$ and $|h|_1$.
Therefore,
\begin{eqnarray*} \int_{-2\pi m''}^{2\pi m''} \bigl|\psi^{(4)}_{\Delta}(u)\bigr|^2\,du
&\leq& C \Delta^2 \biggl[\int x^4 h^2(x) \,dx + \Delta^2\biggl(\int x^2 h^2(x)\,dx + 4
\pi m''|h|_1^{4}\biggr) \\
&&\hspace*{145pt}{} +\Delta^4 \|h\|^2 + 4\pi m'' \Delta^6 \biggr]\\
&\leq& C_1 \Delta^2\biggl[\int x^4h^2(x)\,dx + \Delta^2\int x^2h^2(x)\,dx +\Delta
^4 \|h\|^2\biggr]\\
&&{} + C_2 m'' \Delta^4.
\end{eqnarray*}
Thus, using Assumptions (H1), (H3), (H4),
\[
\int_{[-2\pi m'', 2\pi m'']} \bigl|\psi^{(4)}_{\Delta}(u)\bigr|^2\,du \leq K(\Delta
^2 + m''\Delta^4).
\]
As $m''\Delta^4\leq n\Delta^5$ and $n\Delta^3\leq1$ we get
$\int_{[-2\pi m'', 2\pi m'']} |\psi^{(4)}_{\Delta}(u)|^2\,du\leq 2
K\Delta^2.$
This together with (\ref{presquev}) yields $v=c\sqrt{ m''}/\Delta$
where $c$ is a constant.

Applying Lemma \ref{Concent} yields, for $\epsilon^2=1/2$ and
$p(m,m')$ given by (\ref{p3mm}) yields
\begin{eqnarray*}
&&{\mathbb E} \Bigl(\sup_{t\in
S_m+S_{m'}, \|t\|=1} [\nu_{n,1}^{P}(t)]^2
-p(m,m') \Bigr)_+ \\
&&\qquad \leq C_1 \biggl(\frac{\sqrt{m''}}{n\Delta}
e^{-C_2\sqrt{m''}} +
\frac{k_n^2m''}{n^2\Delta}e^{-C_3\sqrt{n}/k_n} \biggr)
\end{eqnarray*}
as
$p(m,m')=4H^2$. We choose
\[
k_n=\frac{C_3}4
\frac{\sqrt{n}}{\ln(n\Delta)},
\]
and as $m\leq n\Delta$, we get
\begin{eqnarray*}
&&{\mathbb E} \Bigl(\sup_{t\in
S_m+S_{m'}, \|t\|=1} [\nu_{n,1}^{P}(t)]^2
-p(m,m') \Bigr)_+ \\
&&\qquad \leq C'_1 \biggl(\frac{\sqrt{m''}}{n\Delta}
e^{-C_2\sqrt{m''}} + \frac1{(\Delta n)^4\ln^2(n\Delta)} \biggr).
\end{eqnarray*}
%
As $C_2xe^{-C_2x}$ is decreasing for $x\geq1/C_2$, and its maximum is
$1/(eC_2)$, we get
\begin{eqnarray*} \sum_{m'=1}^{m_n} \sqrt{m''}e^{-C_2\sqrt{m''}}&\leq&
\sum_{\sqrt{m'}\leq1/C_2} (eC_2)^{-1} +
\sum_{\sqrt{m'}\geq1/C_2} \sqrt{m'}e^{-C_2\sqrt{m'}}\\ &\leq& \frac
1{eC_2^3} +
\sum_{m'=1}^{\infty} \sqrt{m'}e^{-C_2\sqrt{m'}} <+\infty.
\end{eqnarray*}
It follows that
\[
\sum_{m'=1}^{m_n} {\mathbb E} \Bigl(\sup_{t\in
S_m+S_{m'}, \|t\|=1} [\nu_{n,1}^{P}(t)]^2
-p(m,m') \Bigr)_+ \leq\frac C{n\Delta}.
\]

Let us now study the second term $\nu_{n,j}^{(R)}(t)$ in the
decomposition of $\nu_{n,j}(t)$. The cases $j=3,4$ being similar, we
consider only $\nu_{n,j}^{(R)}(t)$
for $j=1$:
\begin{eqnarray*} && {\mathbb E} \Bigl[\sup_{t\in S_{m_n}, \|t\|=1} \bigl|\nu
_{n,1}^{(R)}(t)\bigr|^2 \Bigr]\\
&&\qquad \leq
\frac1{4\pi^2\Delta^2} {\mathbb E} \Biggl(\int_{-\pi m_n}^{\pi m_n} \biggl|\frac
1n \sum_{k=1}^n
\bigl(Z_k^2\1_{Z_k^2>k_n\sqrt{\Delta}}e^{iuZ_k}-{\mathbb E}(Z_k^2\1
_{Z_k^2>k_n\sqrt{\Delta}}e^{iuZ_k})\bigr)
\biggr|^2\\
&&\hspace*{276pt}
{}\times |\psi_{\Delta}^2(u)|^2\,du \Biggr)\\
&&\qquad \leq \frac{{\mathbb E}(Z_1^4\1_{Z_1^2>k_n\sqrt{\Delta}}) }{4n\pi
^2\Delta^2}
\int_{-\pi m_n}^{\pi m_n} \,du \leq\frac{m_n {\mathbb E}(Z_1^{4+
2p})}{2\pi n\Delta^2(k_n\sqrt{\Delta
})^p }\\
&&\qquad \leq K\frac{{\mathbb E}(Z_1^{4+ 2p}/\Delta)\ln^p(n\Delta)}{2\pi
(n\Delta)^{p/2} },
\end{eqnarray*}
using $m_n\leq n\Delta$ and recalling that $k_n=(C_3/4)(\sqrt{n}/\ln
(n\Delta))$. Taking $p=2$, which
is possible because ${\mathbb E}(Z_1^8)<+\infty$, gives a bound of
order $\ln^2(n\Delta)/(n\Delta)$.

Proposition \ref{concentP} is proved.
\end{pf*}

\begin{pf*}{Proof of Proposition \ref{concentP4}}
For $\nu_{n,2}$,
the variables are bounded without splitting, and the function $f_t$ is
replaced by $\tilde f_t (z)=(2\pi\Delta)^{-1} \langle t^*, e^{iz\cdot}\psi
_{\Delta}''\rangle$. We just check the orders of $M$, $H^2$ and $v$ for
the application of Lemma \ref{Concent}.
For $t\in S_{m''}=S_m+S_{m'}$ and $\|t\|\leq1$, we have
\begin{eqnarray*} \sup_{z\in{\mathbb R}}|\tilde f_t(z)|&\leq& \frac
1{2\pi\Delta}\sqrt{\int_{-\pi m''}^{\pi m''} |t^*(-u)|^2\,du \int_{-\pi
m''}^{\pi m''} |\psi''_{\Delta}(u)|^2\,du} \\
&\leq& \sqrt{m''}\frac{{\mathbb E}(Z_1^2)}{\Delta}\leq C\sqrt{m''}:=M.
\end{eqnarray*}
Next,
\begin{eqnarray*}
{\mathbb E}\Bigl(\sup_{t\in S_m+S_{m'}, \|t\|=1} [\nu_{n,2}(t)]^2\Bigr)
&\leq&
\frac1{2\pi n\Delta^2} \int_{-\pi m''}^{\pi m''}|\psi_{\Delta}''(u)|^2
\,du \\
&\leq&\frac{m'' {\mathbb E}^2(Z_1^2/\Delta)}{n\Delta}:=H^2.
\end{eqnarray*}
Following the same line as previously for $v$, we get
\begin{eqnarray*}
&&\operatorname{Var}(\tilde f_t(Z_1)) \\
&&\qquad \leq \frac1{4\pi^2\Delta^2} \biggl(\iint_{[-\pi
m'', \pi m'']^2} |\psi_{\Delta}(u-v)|^2|\psi_{\Delta}''(u)|^2|\psi
_{\Delta}''(-v)|^2 \,du\,dv \biggr)^{1/2}.
\end{eqnarray*}
As $\psi''_{\Delta}=\Delta(\phi'+\Delta\phi^2]\psi_{\Delta}$, we get
(recall that $M_{\phi}= |b|+|h^*|_1$ is the upper bound of $|\phi(u)|$)
\begin{eqnarray*}
\operatorname{Var}(\tilde f_t(Z_1)) &\leq& \frac1{4\pi
^2\Delta^2}
\int_{-\pi m''}^{\pi m''} |\psi_{\Delta}''(x)|^2\,dx \leq \frac{2\Delta
^2(\|h^*\|^2 + 2\pi m'' \Delta^2M_\phi^2)}{4\pi^2\Delta^2} \\ & \leq&
\frac1\pi(\|h\|^2+ M_{\phi}^2 m_n\Delta^2)
\leq \frac{\|h\|^2+ M_{\phi}^2}\pi:=v
\end{eqnarray*}
as $m_n\Delta^2\leq n\Delta^3\leq1$.
\end{pf*}

\begin{pf*}{Proof of Proposition \ref{concentP12}}
Here, $f_t$ is replaced by $\breve f_t(z)=\break z\1_{|z|\leq k'_n\sqrt{\Delta}}\langle
t^*,
e^{iz\cdot}\psi_{\Delta}'\rangle$.
Using now that $|\psi'_{\Delta}(u)|\leq{\mathbb E}(|Z_1|)\leq\sqrt
{{\mathbb E}(Z_1^2)}$, we obtain here that $M=k'_n\sqrt{m''}\sqrt
{{\mathbb E}(Z_1^2/\Delta)}$. On the other hand, we find
$H^2=m''{\mathbb E}^2(Z_1^2)/(n\Delta^2)$. Last, we find
\begin{eqnarray*}
&&\operatorname{Var}(\breve f_t(Z_1)) \\
&&\qquad \leq \frac1{4\pi^2\Delta^2} \biggl(\iint_{[-\pi
m'', \pi m'']^2} \bigl|\psi^{(2)}_{\Delta}(u-v)\bigr|^2|\psi_{\Delta}'(u)|^2|\psi
_{\Delta}'(-v)|^2 \,du\,dv \biggr)^{1/2}.
\end{eqnarray*}
With the bounds for $|\psi_{\Delta}'|$ and $\int_{-2\pi m''}^{2\pi
m''}|\psi_{\Delta}''(z)|^2\,dz$, we obtain $v= c{\mathbb E}(Z_1^2/\break\Delta
)\sqrt{m''}$.
\end{pf*}

\subsection{\texorpdfstring{Proof of Proposition \protect\ref{rateofp}}{Proof of Proposition 4.2}}
Let us take $m=O((n\Delta)^{1/(2a+1)})$. When $p\in
{\mathcal C}(a,L)$, the first two terms of (\ref{riskcomplet}) are of
order $O((n\Delta)^{-2a/(2a+1)})$. The third term is $O(\Delta^2
m^{2(2-a)_+})$. If $a\geq2$, its order is $\Delta^2$ and is less than
$1/(n\Delta)$ if $n\Delta^3\leq1$.

If $a\in(0,2)$, $\Delta^2m^{2(2-a)} = O(\Delta^2(n\Delta
)^{2(2-a)/(1+2a)})$ which has lower rate than
$O((n\Delta)^{-2a/(2a+1)})$ if $\Delta^2 (n\Delta)^{4/(1+2a)} \leq
O(1)$, that
is
$n\Delta^{1+(1+2a)/2}=n\Delta^{3/2+a}\leq O(1)$.
We must consider in addition the terms $\Delta^2m^3$ and $\Delta^4m^7$.
As previously, $\Delta^2m^3\leq(n\Delta)^{-2a/(2a+1)}$ if $n\Delta
^{(6a+5)/(2a+3)}\leq O(1)$ that is\break $n\Delta^{5/3}\leq1$ if $a>0$ and
$n\Delta^2$ if $a\geq1/2$. Moreover, $\Delta^4m^7\leq(n\Delta
)^{-2a/(2a+1)}$ if $n\Delta^{(10a+11)/(2a+7)}\leq1$ that is $n\Delta
^{11/7}\leq1$ if $a>0$ and $n\Delta^2\leq1$ if $a\geq1/2$.

\subsection{\texorpdfstring{Proof of Proposition \protect\ref{risk2}}{Proof of Proposition 4.1}}
As previously, $\|\bar p_m - p\|^2= \frac1{2\pi} (\|p^*-p_m^*\|^2 + \|
p_m^*-\bar p_m^*\|^2)$.
The variance of $\bar p_m$ satisfies
\begin{eqnarray*}
{\mathbb E}(\|\bar p_m-p_m\|^2) &=& \frac1{2\pi} {\mathbb E}(\|\bar
p_m^*-p_m^*\|^2)\\
&=&\frac1{2\pi} \int_{-\pi m}^{\pi m} \bigl(\operatorname{Var}(\bar p^*(u))+ |{\mathbb
E}(\bar p^*(u))-p^*(u) |^2 \bigr) \,du,
\end{eqnarray*}
where
\[
\operatorname{Var}(\bar p^*(u))\leq\frac{{\mathbb E}(Z_1^6)}{n\Delta^2}=\frac
{{\mathbb E}(Z_1^6/\Delta)}{n\Delta}.
\]

We have $| h^*(u)|\le|h|_1$. By Lemma \ref{bias1}, $|{\tilde\phi}(u)
|\leq|b| + |u| (|h|_1+\sigma^2)\leq C(1+|u|)$.
Inserting these bounds in (\ref{biaspbar}) implies
\begin{eqnarray} \label{bias*}
|{\mathbb E}(\bar p^*(u))-p^*(u)| &\leq&
C\Delta|p^*(u)||u|(1+|u|) \nonumber\\[-8pt]\\[-8pt]
&&{}+ C' \Delta(1 + |u|) +C'' \Delta^2 (1+|u|)^3.\nonumber
\end{eqnarray}
%
Gathering the terms gives the announced bound for the risk of $\bar p_m$.
This ends the proof of Proposition \ref{risk2}.

\subsection{\texorpdfstring{Proof of Theorem \protect\ref{main2}}{Proof of Theorem 4.1}}
The proof follows the same lines as for the adaptive estimator of
$h$. We introduce, for $0<\varrho<1$,
\[
\Omega_b:= \biggl\{ \biggl|
\frac{[(1/(n\Delta))\sum_{k=1}^{n} Z_k^6]}{({\mathbb
E}(Z_1^6/\Delta))} -1 \biggr|\leq\varrho\biggr\}.
\]
Provided that
${\mathbb E}(Z_1^{24})<\infty$, we can make use of the Rosenthal
inequality to obtain:
\[
{\mathbb E}(\|\bar p_{\bar m}-p\|^2\1_{\Omega_\varrho^c})\le C/n\Delta.
\]
For the study of ${\mathbb E}(\|\bar p_{\bar
m}-p\|^2\1_{\Omega_\varrho})$, the decomposition is similar to the
previous case [see (\ref{risk})] where $\hat h_{\hat m}, h$ are
now replaced by $\bar p_{\bar m},p$. The processes $R_n(t)$ and
$\nu_n(t)$ are given by
\[
\nu_n(t)=\frac{1}{2\pi} \langle\bar
p^*-{\mathbb E}(\bar p^*), t^*\rangle,\qquad  R_n(t)=\frac{1}{2\pi} \langle
{\mathbb E}(\bar p^*) - p^*, t^*\rangle.
\]
The term $R_n(t)$ is dealt using (\ref{bias*}). For the term
containing $\nu_n(t)$, we need apply Lemma \ref{Concent}. So, $\nu_n$ is
split into the sum of a principal and a residual term, respectively
denoted by $\nu_n^{P}$ and $\nu_n^{R}$ with
\begin{eqnarray} \label{nunP}
\nu_{n}^P(t) = \frac{1}{n}
\sum_{k=1}^{n} [f_t(Z_k)-{\mathbb E}(f_t(Z_k))] \nonumber\\[-8pt]\\[-8pt]
\eqntext{\displaystyle \mbox{with }
f_t(z)=\frac1{2\pi\Delta} z^3\1_{|z|^3\leq k_n\sqrt{\Delta}} \langle
t^*, e^{iz\cdot}\rangle,}
\end{eqnarray}
and $\nu_{n}^R(t)=\nu_{n}(t) - \nu_{n}^P(t)$. Everything is
analogous. The difference is that, for applying Lemma \ref{Concent}, we
have to bound $\int_{-2\pi m''}^{2\pi m''} |\psi_\Delta^{(6)}(u)|^2
\,du$ (instead of $\int_{-2\pi m''}^{2\pi m''} |\psi_\Delta^{(4)}(u)|^2
\,du$ previously). Using $\psi_{\Delta}'=\Delta\tilde\phi
\psi_{\Delta}$ [see (\ref{phidu})--(\ref{phitilde})], we find
\begin{eqnarray*}
\psi_\Delta^{(6)}&=& \Delta \psi_\Delta\phi^{(5)} + \Delta^2 \psi
_\Delta\bigl[6
\tilde\phi\phi^{(4)} +
15\phi^{(3)} \bigl(\phi'(u)-\sigma^2\bigr)\bigr]\\
&&{}+\Delta^3 \psi_\Delta \bigl[15
\phi^{(3)}\tilde\phi^2 +
60 \phi'' \bigl(\phi'(u)-\sigma^2\bigr) \tilde\phi+ 15 \bigl(\phi'(u)-\sigma^2\bigr)^3\bigr]\\
&& {}+
\Delta^4 \psi_\Delta\bigl[17 \phi''{\tilde\phi}^{(3)}+ 36 \tilde\phi^{(2)}
\bigl(\phi'(u)-\sigma^2\bigr)^{2} \bigr]\\
&&{} + 12\Delta^5 \psi_\Delta\tilde\phi^4 \bigl(\phi'(u)-\sigma^2\bigr) + \Delta
^6 \psi_\Delta\tilde\phi^6.
\end{eqnarray*}
%
Now, $\tilde\phi(u)\leq C(1+|u|)$ and all the derivatives of $\tilde\phi
, \phi$ are bounded. Moreover, under (H6), $\int|\phi^{(5)}(u)|^2\,du =
\int x^6|p(x)|^2 \,dx <+\infty$. Thus, we find
the following bound:
\[
\int_{-2\pi m''}^{2\pi m''}|\psi^{(6)}_\Delta|^2 \leq C\Delta^2( 1 +
\Delta^2 m^3+ \Delta^4 m^5+ \Delta^6 m^7+\Delta^8 m^9 + \Delta^{10}
m^{13})=O(\Delta^2),
\]
as $m\leq\sqrt{n\Delta}$. The proof may then be completed as for $\hat
h_{\hat m}$.

\subsection{\texorpdfstring{Proof of Proposition \protect\ref{gamma}}{Proof of Proposition 5.1}}
Proof of (i). The assumptions and the fact that $r\le1$ imply
\[
|\Gamma_{\Delta}|^r=\biggl|\sum_{s\le\Delta}\Gamma_s-\Gamma_{s_{-}}\biggr|^r\le
\sum_{s\le\Delta}|\Gamma_s-\Gamma_{s_{-}}|^r.
\]
Taking expectations yields ${\mathbb E}|\Gamma_{\Delta}|^r\le
\Delta\int|\gamma|^r n(\gamma)\,d\gamma.$

 Proof of (ii). Consider $f$ a nonnegative function such
that $f(0)=0$. We have
\[
{\mathbb E}\sum_{s\le t}f(X_s - X_{s_-})= {\mathbb E}\sum_{s\le
t}f(B_{\Gamma_s} -
B_{\Gamma_{s_-}}).
\]
Then,
$
\sum_{s\le t} {\mathbb E}f(B_{\Gamma_s} -
B_{\Gamma_{s_-}})=\sum_{s\le t} \int_{{\mathbb R}} f(x) ({\mathbb E}
e^{(-x^2/2(\Gamma_s - \Gamma_{s_-}))}
\frac{1}{\sqrt{2 \pi(\Gamma_s - \Gamma_{s_-})}} )\,dx.
$
Since, for all $x$,
\[
{\mathbb E}\sum_{s\le t} e^{(-x^2/2(\Gamma_s - \Gamma_{s_-}))}
\frac{1}{\sqrt{2 \pi(\Gamma_s - \Gamma_{s_-})}} = t \int_0^{+\infty}
e^{-x^2/2\gamma}
\frac{1}{\sqrt{2 \pi\gamma}}n_\Gamma(\gamma)\,d\gamma,
\]
we get the formula for $n_X$. Setting $m_{\alpha}={\mathbb E}|X|^{\alpha
}$, for $X$ a standard Gaussian variable, yields
\[
\int_{\mathbb R}|x|^{\alpha} n_X(x) \,dx= m_{\alpha} \int_0^{+\infty}
\gamma^{\alpha/2}n_{\Gamma}(\gamma)\,d\gamma.
\]
Thus,
$
{\mathbb E}|X_{\Delta}|^r= m_r{\mathbb E}(\Gamma_{\Delta}^{r/2}).
$
As $r/2 \le1$,
$
\Gamma_{\Delta}^{r/2}=(\sum_{s\le\Delta}\Gamma_s-\Gamma
_{s_{-}})^{r/2}\le\break\sum_{s\le\Delta}(\Gamma_s-\Gamma_{s_{-}})^{r/2}.
$
Taking expectation gives the result.

 Proof of (iii). The result is proved, for example, in
Barndorff-Nielsen, Shephard and Winkel [(\citeyear{BSW2006}), Theorem 1,
page 804] [see also A\"{i}t-Sahalia and Jacod (\citeyear{AJ2007})].

\subsection{\texorpdfstring{Proof of Proposition \protect\ref{emp}}{Proof of Proposition 5.2}}
We have ${\mathbb E}(Z_k)=\Delta b$ and, for $\ell\ge2$,
${\mathbb E}(Z_k^{\ell})=\Delta c_{\ell}+ o(\Delta)$. Therefore, $\hat
b$ is an unbiased estimator of $b$ and, for $\ell\ge2$, $\sqrt{n\Delta
}|{\mathbb E}{\hat c}_{\ell}-c_{\ell}| = \sqrt{n\Delta}O(\Delta)$.
Hence, the additional condition $n\Delta^{3}=o(1)$ to erase the bias.

Setting $c_1=b$, $\hat c_1=\hat b$, as
$\operatorname{Var}Z_k^{\ell}= \Delta c_{2\ell}+ o(\Delta)$ for $\ell\geq1$, we have
$
n\Delta\operatorname{Var}{\hat c}_{\ell}= c_{2\ell}+ O(\Delta).
$
Writing
$
\sqrt{n\Delta}({\hat c}_{\ell}-{\mathbb E}{\hat c}_{\ell})= (n\Delta
)^{-1/2}\sum_{k=1}^{n} (Z_k^{\ell}- {\mathbb E}Z_k^{\ell})= \sum
_{k=1}^{n} \chi_{k,n},
$
it is now enough to prove that $\sum_{k=1}^{n}{\mathbb E}| \chi
_{k,n}|^{2+\varepsilon}$ tends to $0$.
Under the assumption, we have
\[
\sum_{k=1}^{n}{\mathbb E}|\chi_{k,n}|^{2+\varepsilon}\le \frac
{C}{n^{\varepsilon/2}\Delta^{1+\varepsilon/2}} \bigl({\mathbb E}|Z_k|^{\ell
(2+\varepsilon)}+ |{\mathbb E}(Z_k^{\ell})|^{2+\varepsilon} \bigr) \le\frac
{C}{(n \Delta)^{\varepsilon/2}},
\]
which gives the result.

\subsection{\texorpdfstring{Proof of Proposition \protect\ref{estimsigma}}{Proof of Proposition 5.3}}
The study of (\ref{sigmahatnr}) relies on the following result which is
standard for $r=2$.
\begin{lem} Let $Y_t= \theta t + \sigma W_t$ for $\theta$ a constant
and consider
$
{\tilde\sigma}_n^{(r)}= \frac{1}{m_r n \Delta^{r/2}} \sum_{k=1}^{n}
|Y_{k\Delta}-Y_{(k-1)\Delta}|^r.
$

Then, for all $r$, $\sqrt{n}({\tilde\sigma}_n^{(r)}-\sigma^{r})$
converges in distribution to a
centered Gaussian distribution with variance $ \sigma
^{2r}(m_{2r}/m_r^{2} -1)$ as $n$ tends to infinity, $\Delta$ tends to
$0$, $n\Delta$ tends to infinity, and $n\Delta^2$ tends to $0$.
\end{lem}

\begin{pf} We have
$
{\mathbb E}{\tilde\sigma}_n^{(r)}= \frac{1}{m_r}{\mathbb E}|\theta\sqrt
{\Delta}+\sigma X|^r,
$
for $X$ a standard Gaussian variable. Thus,
%
\begin{eqnarray*}
{\mathbb E}{\tilde\sigma}_n^{(r)}-\sigma^r&=&\sigma^r (e^{-\theta^2\Delta
/2\sigma^2}-1 )\\
&&{}+\frac{1}{m_r}e^{-\theta^2\Delta/2\sigma^2}\int|u|^{r}
(e^{\theta u \sqrt{\Delta}/\sigma^2}-1)e^{-{u^2}/{(2\sigma^2)}}\frac
{du}{\sigma\sqrt{2\pi}}.
\end{eqnarray*}
Noting that
$
e^{\theta u \sqrt{\Delta}/\sigma^2}-1=\theta u \sqrt{\Delta}/\sigma^2 +
\Delta\sum_{n\ge2} \frac{1}{n!}(u \theta/\sigma^2)^n \Delta^{n/2 -1}
$
and that $\int|u|^r u e^{-{u^2}/{(2\sigma^2)}}\,du/(\sigma\sqrt{2\pi
})=0$, we easily obtain
\[
\bigl|{\mathbb E}{\tilde\sigma}_n^{(r)}-\sigma^r\bigr| \le c\Delta.
\]
Thus, $\sqrt{n}|{\mathbb E}{\tilde\sigma}_n^{(r)}-\sigma^r|=o(1)$ if
$\sqrt{n}\Delta= (n\Delta^2)^{1/2}=o(1)$.
Noting that ${\mathbb E}|\theta\sqrt{\Delta}+\sigma X|^k$ converges to
$\sigma^k m_k$ as $\Delta$ tends to $0$, we get
$
n \operatorname{Var}{\tilde\sigma}_n^{(r)} \rightarrow\sigma
^{2r}(m_{2r}/\break m_r^2 -1).
$

Finally, we look at
$
\chi_{k,n}= n^{-1} ( |\theta\sqrt{\Delta}+ \sigma(W_{k\Delta}
-W_{(k-1)\Delta})/\sqrt{\Delta}|^r - {\mathbb E}|\theta\sqrt{\Delta
}+\sigma X|^r ),
$
which satisfies
$
n {\mathbb E}\chi_{k,n}^{4} \le c/n^3.$
Hence, $\sqrt{n}({\tilde\sigma}_n^{(r)} - {\mathbb E} {\tilde\sigma
}_n^{(r)} )$ converges in distribution to the centered Gaussian with
the announced variance which completes the proof.
\end{pf}

\begin{pf*}{Proof of \textup{(i)}}
As noted above, $L_t= b_0 t + \sigma W_t +
\Gamma_t$ with $b_0=b-\int xn(x)\,dx$. Using that, for $r \le1$, $||\sum
a_i+b_i|^r - |\sum a_i|^r|\le\sum|b_i|^r$, we get
$
|{\hat\sigma}_n^{(r)}- {\tilde\sigma}_n^{(r)}| \le\frac{1}{m_r n
\Delta^{r/2}} \sum_{k=1}^{n} |\Gamma_{k\Delta}-\Gamma_{(k-1)\Delta}|^r,
$
where ${\tilde\sigma}_n^{(r)}$ is built with $Y_t=b_0 t +\sigma
W_t$ as in the previous lemma. Thus, applying Proposition \ref{gamma}(i),
\[
{\mathbb E}\sqrt{n}\bigl|{\hat\sigma}_n^{(r)}- {\tilde\sigma}_n^{(r)}\bigr| \le
\frac{1}{m_r }\sqrt{n}\Delta^{1-r/2} \int|x|^{r}n(x)\,dx.
\]
Since $r <1$, the constraint $n\Delta^{2-r}=o(1)$ can be fulfilled and
implies $n\Delta^2=o(1)$. Hence, the result follows from the previous
proposition.

\textsc{Proof of} (ii). The proof is analogous to the previous one
[using Proposition~\ref{gamma}(ii)] and is omitted. As $\sigma
(r)=[\hat\sigma_n^{(r)}]^{1/r}$, we conclude for $\hat\sigma(r)$ by
using the delta-method.
\end{pf*}

\begin{appendix}
\section*{Appendix: The Talagrand inequality} \label{Ap}
\setcounter{lem}{0}
The following result
follows from the Talagrand concentration inequality given in Klein
and Rio (\citeyear{KR2005}) and arguments in Birg\'{e} and Massart
(\citeyear{BM1998}) (see the
proof of their Corollary 2, page 354).
\begin{lem}[(Talagrand inequality)]
\label{Concent}  Let $Y_1, \dots, Y_n$ be
independent random
variables, let $\nu_{n,Y}(f)=(1/n)\sum_{i=1}^n [f(Y_i)-{\mathbb
E}(f(Y_i))]$ and let ${\mathcal F}$ be a countable class of
uniformly bounded measurable functions. Then for $\epsilon^2>0$
\begin{eqnarray*}
&& \mathbb{E} \biggl[\sup_{f\in{\mathcal
F}}|\nu_{n,Y}(f)|^2-2(1+2\epsilon^2)H^2 \biggr]_+\\
&&\qquad \leq\frac
4{K_1} \biggl(\frac vn e^{-K_1\epsilon^2 nH^2/v} +
\frac{98M^2}{K_1n^2C^2(\epsilon^2)} e^{-2K_1
C(\epsilon^2)\epsilon/{(7\sqrt{2})}{nH}/{M}} \biggr),
\end{eqnarray*}
with $C(\epsilon^2)=\sqrt{1+\epsilon^2}-1$, $K_1=1/6$  and
\[
\sup_{f\in{\mathcal F}}\|f\|_{\infty}\leq M,\qquad
\mathbb{E} \Bigl[\sup_{f\in{\mathcal F}}|\nu_{n,Y}(f)| \Bigr]\leq H,\qquad
\sup_{f\in{\mathcal F}}\frac{1}{n}\sum_{k=1}^n\operatorname{Var}(f(Y_k)) \leq v.
\]
\end{lem}

By standard density arguments, this result can be extended to the case
where ${\mathcal F}$ is a unit ball of a linear normed space, after
checking that $f\mapsto\nu_n(f)$ is continuous and ${\mathcal F}$
contains a countable dense family.
\end{appendix}

\printaddresses

\end{document}